\documentclass[twoside,centertags]{amsart}
\usepackage[cmtip,color,all]{xy}
\usepackage{amsmath}
\usepackage{amsfonts}
\usepackage{amssymb}
\usepackage{amsthm}
\usepackage{graphicx}
\usepackage{hyperref}
\usepackage{color}
\usepackage{mathrsfs}
\usepackage[shortlabels]{enumitem}
\usepackage[mathscr]{euscript}

\newtheorem{theorem}[equation]{Theorem}
\newtheorem{corollary}[equation]{Corollary}
\newtheorem{lemma}[equation]{Lemma}

\newtheorem{proposition}[equation]{Proposition}
\theoremstyle{definition}
\newtheorem{definition}[equation]{Definition}
\newtheorem{example}[equation]{Example}

\newtheorem{remark}[equation]{Remark}
\theoremstyle{remark}

\newcommand{\s}{\textsf{s}}
\def\facto{\mathscr{E}}

\newcommand{\Ob}{\operatorname{Ob}}

\newcommand{\C}{\mathscr{C}}
\newcommand{\D}{\mathscr{D}}
\newcommand{\Fun}{\mathrm{Fun}}
\newcommand{\h}{\mathrm{h}}
\newcommand{\Ar}{\mathrm{Ar}}
\renewcommand{\S}{\mathrm{S}}
\newcommand{\catinf}{\mathsf{Cat_{\infty}}}
\newcommand{\catn}{\mathsf{Cat_{n}}}
\newcommand{\dia}{\mathsf{Dia}}
\renewcommand{\c}{\mathsf{C}}

\newcommand{\x}{\mathsf{X}}
\newcommand{\y}{\mathsf{Y}}
\newcommand{\w}{\mathsf{W}}
\newcommand{\T}{\mathcal{T}}

\newcommand{\z}{\mathsf{Z}}

\renewcommand{\a}{\mathsf{A}}
\newcommand{\can}{\mathrm{can}}
\newcommand{\diag}{\mathrm{dia}}
\def\colim{\operatorname{colim}}
\newcommand{\dirf}{\mathcal{D}ir_{f}}
\newcommand{\PreDer}{\mathsf{PreDer}}
\newcommand{\Der}{\mathsf{Der}}

\def\ob{\operatorname{Ob}}

\def\Map{\mathrm{Map}}
\def\der{\mathbb{D}}

\def\op{\operatorname{op}}

\def\sk{\mathrm{sk}}

\numberwithin{equation}{section}

\begin{document}
\title[Higher homotopy categories and $K$-theory]{Higher homotopy categories, \\ higher derivators, and $K$-theory} %
\author[G. Raptis]{George Raptis}
\address{\newline
G. Raptis \newline
Fakult\"{a}t f\"ur Mathematik \\
Universit\"{a}t Regensburg \\
93040 Regensburg, Germany}
\email{georgios.raptis@ur.de}

\begin{abstract}
For every $\infty$-category $\C$, there is a homotopy $n$-category $\h_n \C$ and a canonical functor $\gamma_n \colon \C \to \h_n \C$. We study these higher homotopy categories, 
especially in connection with the existence and preservation of (co)limits, by introducing a higher categorical notion of weak colimit. Using homotopy $n$-categories, we 
introduce the notion of an $n$-derivator and study the main examples arising from $\infty$-categories. Following the work of Maltsiniotis and Garkusha, we define $K$-theory for 
$\infty$-derivators and prove that the canonical comparison map from the Waldhausen $K$-theory of $\C$ to the $K$-theory of the associated $n$-derivator $\der_{\C}^{(n)}$ is 
$(n+1)$-connected. We also prove that this comparison map identifies derivator $K$-theory of $\infty$-derivators in terms of a universal property. Moreover,  using the canonical structure of higher weak pushouts in the homotopy $n$-category, we also define a $K$-theory space $K(\h_n \C, \can)$ associated to $\h_n \C$. We prove that the canonical comparison map from the Waldhausen $K$-theory of $\C$ to $K(\h_n \C, \can)$ is $n$-connected.
\end{abstract}

\maketitle
\setcounter{tocdepth}{1}
\tableofcontents

\section{Introduction}

It has been long understood in homotopy theory that the homotopy category is only a crude invariant of a much richer homotopy-theoretic structure. The problem of finding a suitable formalism 
for this additional structure, one that encodes homotopy-theoretic extensions of ordinary categorical notions, led to several foundational approaches, each with its own distinctive advantages
and special  characteristics. The theory of $\infty$-categories (or quasi--categories) \cite{Joyal_1, LurHA, LurHTT, Cisinski} is one of several successful approaches to develop useful foundations for the study of homotopy theories and has led to groundbreaking new perspectives and results in the field. 

\smallskip

Even though passing to the homotopy category certainly neglects homotopy-theoretic information, the general problem of understanding how much information this process retains still poses interesting questions in specific contexts. This has inspired many important developments, for example, in the context of rigidity theorems for homotopy theories \cite{Franke, Schwede2, Roitzheim, Patchkoria},
derived/homotopical Morita and tilting theory (see \cite{Schwede} for a nice survey), or in connection with $K$-theory regarded as an invariant of homotopy theories \cite{DS, Neeman, Schlichting, Malt}. 

\smallskip

The theory of derivators, first introduced and developed by Grothendieck \cite{derivateurs}, Heller \cite{Heller} and Franke \cite{Franke}, is a different foundational approach, based on the idea of considering the homotopy categories of all diagram categories as a remedy to the 
defects of the homotopy category (see also \cite{Groth}). By supplementing the homotopy category with the network of all these (homotopy) categories, it is possible to encode the collection of homotopy (co)limit functors and general homotopy Kan extensions as an enhancement of the homotopy category. 
This approach provides a different, lower (= 2--) categorical formalism for expressing homotopy-theoretic universal properties (see \cite{Cisinski2} and \cite{Tabuada, C-T} for some interesting applications). Moreover, Maltsiniotis \cite{Malt} defined $K$-theory in the 
context of derivators with a view towards partially recovering Waldhausen $K$-theory from the derivator. The $K$-theory of derivators and its 
comparison with Waldhausen $K$-theory has been studied extensively in \cite{C-N, Garkusha2, Garkusha1, Malt, Muro, MR2, MR1}. In the context of the theory of derivators, the question about 
the information retained by the homotopy category is then upgraded to the analogous question for the derivator. The theory of derivators still does not provide in general faithful representations of homotopy theories, however, it is possible in certain cases to recover in a non-canonical way the homotopy theory from the derivator  (see \cite{Renaudin}). 

\smallskip

The purpose of this paper is to extend these ideas on the comparison between homotopy theories and homotopy categories or derivators to $n$-categories (= $(n,1)$-categories),
where the ordinary homotopy category is now replaced by the homotopy $n$-category of an $\infty$-category. More specifically:
\begin{itemize}
 \item[(a)] \emph{Higher homotopy categories.} Using the definition of the higher homotopy categories by Lurie \cite{LurHTT}, we consider the tower of homotopy $n$-categories $\{\h_n \C\}_{n \geq 1}$ associated to an $\infty$-category $\C$, together with the canonical (localization) functors $\gamma_n \colon \C \to \h_n\C$, and we analyse the properties of $\h_n\C$ inherited from $\C$.
 (Sections \ref{sec:higher_homotopy_cats}--\ref{sec:higher_weak_colimits}, \ref{stable-n-categories}--\ref{stable-n-categories2})
 \item[(b)] \emph{Higher derivators.} We introduce a higher categorical notion of a derivator which takes values in $n$-categories. Then we develop the basic theory of higher derivators with a special emphasis on the examples which arise from $\infty$-categories. (Section \ref{sec:higher_derivators})
 \item[(c)] \emph{K-theory of higher derivators.} We extend the definition of derivator $K$-theory by Maltsiniotis \cite{Malt} and Garkusha \cite{Garkusha2, Garkusha1} to $n$-derivators and 
 study the comparison map from Waldhausen $K$-theory. Our main result shows that the comparison map is $(n+1)$-connected (Theorem \ref{comparison_K-theory}). Moreover, following \cite{MR2}, we prove that this comparison map 
 has a universal property (Theorem \ref{universal_prop_derivator_K-theory}). 
 \item[(d)] \emph{K-theory of homotopy n-categories.} In analogy with the $K$-theory of triangulated categories \cite{Neeman}, we introduce $K$-theory for $n$-categories equipped with 
 distinguished squares. In the case of a homotopy $n$-category, we study the comparison map from Waldhausen $K$-theory and prove that it is $n$-connected (Theorem \ref{comparison_K-theory_ncat}). 
\end{itemize}

\smallskip

\noindent (a) \textbf{Higher homotopy categories.} Every $\infty$-category $\C$ has an associated homotopy $n$-category $\h_n \C$ and a canonical functor 
$\gamma_n \colon \C \to \h_n \C$. The construction of the homotopy $n$-category and its properties are studied in \cite{LurHTT}. We will review this construction and its properties 
in Section \ref{sec:higher_homotopy_cats}. Intuitively, for $n \geq 1$, $\h_n \C$ is an $\infty$-category with the same objects as $\C$ and whose mapping spaces are the appropriate Postnikov truncations of the mapping spaces in $\C$. For $n=1$, the homotopy category $\h_1 \C$ is the ordinary homotopy category of $\C$. The collection of homotopy $n$-categories defines a tower of $\infty$-categories:
$$
\xymatrix{ 
& \C \ar[d]_{\gamma_n} \ar[dr]_{\gamma_{n-1}} \ar[drrr]^{\gamma_1} &&&  \\
\cdots \ar[r] & \h_n \C \ar[r] & \h_{n-1} \C \ar[r] & \cdots \ar[r] & \h_{1} \C
}
$$
which approximates $\C$. One of the defects of the homotopy category $\h_1 \C$, which is essentially what the theory of derivators tries to rectify, is that it does not in general inherit 
(co)limits from $\C$. As a general rule, if $\C$ admits (co)limits, then $\h_1 \C$ admits only weak (co)limits -- which may or may not 
be induced from $\C$. (We recall that a weak colimit of a diagram is a cocone on the diagram which admits a morphism to every other such cocone,
but this morphism may not be unique in general.) Do the higher homotopy categories have \emph{better} inheritance properties for (co)limits, and in what sense? This question is closely related 
to the problem of understanding how much information $\h_n \C$ retains from $\C$. We introduce a higher categorical version of weak (co)limit in order to address this question. In analogy with ordinary weak (co)limits, a higher weak colimit is a 
cocone for which the mapping spaces to other cocones are highly connected, but not necessarily contractible. The relative strength of the weak colimit is measured by how highly connected these mapping spaces are; the connectivity of these mapping spaces determines the \emph{order} of the weak colimit. The properties of higher weak (co)limits will be discussed in Section \ref{sec:higher_weak_colimits}. For any simplicial set $K$, there is a canonical functor
$$\Phi_n^K \colon \h_n(\C^K) \to \h_n(\C)^K$$
which is usually not an equivalence. The properties of this functor are relevant for understanding the interaction of $K$-colimits in $\C$ and in $\h_n \C$. One of our conclusions (Corollary \ref{n-d_homotopy_cat}) is the following: 
\smallskip
\begin{center} 
$\Phi_n^K$ induces an equivalence: $\h_{n-\mathrm{dim(}K\mathrm{)}}(\C^K) \simeq \h_{n-\mathrm{dim(}K\mathrm{)}}\big(\h_n(\C)^K\big)$.
\end{center}
\smallskip
Moreover, in connection with higher weak colimits in $\h_n \C$, we also conclude (Corollary \ref{cor_weak_pushouts}):
\smallskip
\begin{center}
Suppose that $\C$ admits finite colimits. Then $\h_n \C$ admits finite coproducts and weak pushouts of order $n-1$. In addition, the functor 
$\gamma_n \colon \C \to \h_n \C$ preserves coproducts and sends pushouts in $\C$ to weak pushouts of order $n-1$. 
\end{center}
\smallskip
These properties of homotopy $n$-categories $\h_n\C$ single out a class of $n$-categories for each $n \geq 1$, called \emph{(finitely) weakly cocomplete 
$n$-categories} (Definition \ref{weakly-cocomplete-n-cat}). These classes of $\infty$-categories form a sequence of refinements between ordinary categories with (finite) coproducts and weak pushouts and (finitely) cocomplete $\infty$-categories. 
We explore further the properties of these $n$-categories in connection with adjoint functor theorems and (higher) Brown representability in joint work with H. K. Nguyen and C. Schrade \cite{more-gaft} (building on and extending our previous joint work in \cite{gaft}). 

\smallskip

\noindent (b) \textbf{Higher derivators.} The main example of a (pre)derivator is given by the 2-functor which sends a small category $I$ to the homotopy category $\h_1(\C^{N(I)})$, 
for suitable choices of an $\infty$-category $\C$. Using homotopy $n$-categories instead, we may consider more generally the example of the enriched functor which sends a simplicial set $K$ to the 
homotopy $n$-category $\h_n(\C^K)$. Following the axiomatic definition of ordinary derivators, we introduce a definition of an $n$-derivator, which is meant to encapsulate the main abstract properties of this example. The basic definitions and properties of (left, right, pointed, stable) $n$-(pre)derivators, $1 \leq n \leq \infty$, will be discussed in Section \ref{sec:higher_derivators}. For any $\infty$-category $\C$, there is an associated $n$-prederivator
$$\der^{(n)}_{\C} \colon \dia^{\op} \to \catn, \ \ K \mapsto \h_n(\C^K),$$
where $\dia$ denotes a category of diagram shapes and $\catn$ is the $\infty$-category of $n$-categories. These assemble to define a tower of $\infty$-prederivators:
$$
\xymatrix{ 
& \der^{(\infty)}_{\C} \ar[d] \ar[dr]  \ar[drrr] &&&  \\
\cdots \ar[r] & \der^{(n)}_{\C} \ar[r] & \der^{(n-1)}_{\C} \ar[r] & \cdots \ar[r] & \der^{(1)}_{\C}
}
$$
which approximates $\der^{(\infty)}_{\C} \colon K \mapsto \C^K$. The $n$-prederivator $\der_{\C}^{(n)}$ is an $n$-derivator if certain homotopy Kan extensions exist in $\C$ 
and the corresponding base change transformations satisfy the Beck--Chevalley condition. For any fixed $n \in \mathbb{Z}_{\geq 1} \cup \{\infty\}$, we prove the following fact which can be used to obtain many examples of $n$-derivators 
(Proposition \ref{represented_preder_prop} and Theorem \ref{deriv_criterion}):
\smallskip
\begin{center}
$\C$ admits limits and colimits indexed by diagrams in $\dia$ \newline if and only if $\der_{\C}^{(n)}$ is an $n$-derivator. 
\end{center}
\smallskip
The motivation for higher derivators is to bridge the gap between $\infty$-categories and derivators by introducing a hierarchy of intermediate notions, a different one for each 
categorical level, starting with ordinary derivators. For any fixed $1 \leq n < \infty$, the theory of $n$-derivators is still not a faithful representation of homotopy 
theories; it dwells in an $(n+1)$-categorical context in the same way the classical theory of derivators is restricted to a 2-categorical context. In this respect, our approach using 
$n$-derivators remains close to the original idea of a derivator, and differs from other recent perspectives on (pre)derivators, in which (pre)derivators are reconsidered and revised into 
a model for the theory of $\infty$-categories \cite{FKKR, Carlson, MR2}. We will address the problem of comparing suitable nice classes of $\infty$-categories with $n$-derivators
in future joint work with D.--C. Cisinski.

\smallskip

\noindent (c) \textbf{K-theory of higher derivators.} $K$-theory for (pointed, right) derivators was introduced by Maltsiniotis \cite{Malt} and Garkusha \cite{Garkusha2, Garkusha1}. 
The basic feature of a derivator that allows this definition of $K$-theory is that there is a natural notion of cocartesian square for a derivator. The motivation for introducing 
this $K$-theory is connected with the problem of recovering Waldhausen $K$-theory from the derivator. Maltsiniotis \cite{Malt} conjectured that derivator $K$-theory of stable derivators satisfies 
\emph{additivity}, \emph{localization}, and that it \emph{agrees} with Quillen $K$-theory for exact categories. Cisinski and Neeman \cite{C-N} proved additivity for the derivator $K$-theory of stable derivators and Coley \cite{Coley} has recently extended this result to the unstable context. In joint work with Muro \cite{MR1}, we proved that localization and agreement with Quillen $K$-theory cannot both hold. On the other hand, Muro 
\cite{Muro} proved that agreement with Waldhausen $K$-theory holds for $K_0$ and $K_1$ (see also \cite[Section 6]{Malt}). Moreover, Garkusha \cite{Garkusha1} obtained further positive results in the case of abelian categories. In Section \ref{sec:K-theory_higher_derivators}, after a short review of the $K$-theory of $\infty$-categories, we will define derivator $K$-theory for general (pointed, right) $\infty$-derivators. For any pointed $\infty$-category $\C$ with finite colimits, there is a comparison map to derivator $K$-theory,
$$\mu_n \colon K(\C) \to K(\der_{\C}^{(n)}),$$
and these comparison maps assemble to give a tower of derivator $K$-theories and comparison maps:
$$
\xymatrix{ 
& K(\C) \ar[d]_{\mu_n} \ar[dr]_{\mu_{n-1}}  \ar[drrr]^{\mu_1} &&&  \\
\cdots \ar[r] & K(\der^{(n)}_{\C}) \ar[r] & K(\der^{(n-1)}_{\C}) \ar[r] & \cdots \ar[r] & K(\der^{(1)}_{\C})
}
$$
which approximates $K(\C)$. Our main result on the comparison map $\mu_n$ is the following connectivity estimate (Theorem \ref{comparison_K-theory}):
\smallskip
\begin{center}
 $\mu_n$ is $(n+1)$-connected.
\end{center}
\smallskip
We believe that this connectivity estimate is best possible in general (Remarks \ref{best_possible?} and \ref{t-structure}). Following the ideas of \cite{MR2}, 
we also consider the Waldhausen $K$-theory $K^{W, \ob}(\der)$ (and $K^W(\der)$) of a general (pointed, right) $\infty$-derivator $\der$. This $K$-theory always agrees with the usual $K$-theory (Proposition \ref{Waldhausen_K_theory_der}); the proof is based on a version of the $\s_{\bullet}$-construction in the $\infty$-categorical context (Proposition \ref{comparison_s_constr}). However, Waldhausen $K$-theory of derivators is not invariant under equivalences of derivators in general. Similarly to the case of ordinary derivators treated in \cite{MR2}, 
we prove that the comparison map to derivator $K$-theory,
$$\mu \colon K^{W, \ob}(\der) \to K(\der),$$
identifies derivator $K$-theory in terms of a universal property (Theorem \ref{universal_prop_derivator_K-theory}):
\smallskip
\begin{center}
$\mu$ is the initial natural transformation to a functor\newline which is invariant under equivalences of $\infty$-derivators. 
\end{center}
\smallskip

\smallskip

\noindent (d) \textbf{K-theory of homotopy $n$-categories.} The motivation for introducing $K$-theory for homotopy $n$-categories is to identify the part of Waldhausen $K$-theory which may be 
reconstructed from the homotopy $n$-category. As a basic instance of this phenomenon, we recall that $K_0(\C)$ can be recovered from the triangulated category 
$\h_1 \C$ for any stable $\infty$-category $\C$. The main feature of the homotopy $n$-category that is needed for our definition of $K$-theory is the collection of higher weak pushouts which 
come from pushouts in $\C$. We will revisit the properties of homotopy $n$-categories 
in Subsections \ref{stable-n-categories}--\ref{stable-n-categories2} and discuss possible axiomatizations of these properties. We will then define $K$-theory for pointed $n$-categories with distinguished squares -- this is a higher categorical, but much more elementary, version of Neeman's $K$-theory of categories with squares \cite{Neeman}. For a pointed $\infty$-category $\C$ with finite colimits, we consider 
the $K$-theory $K(\h_n\C, \can)$ associated to $\h_n\C$ with the canonical structure of higher weak pushouts as distinguished squares. For every $n \geq 1$, there is a comparison map 
$$\rho_n \colon K(\C) \to K(\h_n\C, \can),$$
and these assemble to define a tower of $K$-theories and comparison maps:
$$
\xymatrix{ 
& K(\C) \ar[d]_{\rho_n} \ar[dr]_{\rho_{n-1}}  \ar[drrr]^{\rho_1} &&&  \\
\cdots \ar[r] & K(\h_n\C, \can) \ar[r] & K(\h_{n-1} \C, \can) \ar[r] & \cdots \ar[r] & K(\h_1 \C, \can)
}
$$
which approximates $K(\C)$. Our main result in Section \ref{sec:K-theory_higher_homotopy_cats} on the comparison map $\rho_n$ is the following connectivity estimate (Theorem \ref{comparison_K-theory_ncat}):
\smallskip
\begin{center}
$\rho_n$ is $n$-connected.
\end{center}
\smallskip
This connectivity estimate is best possible in general (Remark \ref{best_possible}). Let $P_n X$ denote the Postnikov $n$-truncation of a topological space $X$, that is,  
the homotopy groups of $P_n X$ vanish in degrees $> n$ and the canonical map $X \to P_n X$ is $(n+1)$-connected. Based on the connectivity estimate above, we conclude 
(Corollary \ref{antieau_conj}):
\begin{center}
$P_{n-1} K(\C)$ depends only on $(\h_n \C, \can)$. 
\end{center}
This confirms a recent conjecture of Antieau \cite[Conjecture 8.35]{Antieau} in the case of connective $K$-theory. 

\smallskip

\noindent \textbf{Acknowledgements.} I would like to thank Benjamin Antieau, Denis--Charles Cisinski, Fernando Muro, and Hoang Kim Nguyen for interesting discussions and for their interest in this work. I also thank Martin Gallauer and Christoph Schrade for their interest and helpful comments. In addition, I thank the anonymous referee for their careful reading and useful comments. This work was partially supported by the \emph{SFB 1085 -- Higher Invariants} (University of Regensburg) funded by the DFG. 

\section{Higher homotopy categories} \label{sec:higher_homotopy_cats}

\subsection{$n$-categories} \label{2.1} We recall the definition and basic properties of $n$-categories following \cite[2.3.4]{LurHTT}. Let $\C$ be an $\infty$-category and let $n \geq -1$ be an integer. $\C$ is an $n$-\emph{category} if it satisfies the
following conditions:
\begin{itemize}
 \item[(1)] Given a pair of maps $f,f': \Delta^n \to \C$, if $f$ and $f'$ are homotopic relative to $\partial \Delta^n$, then $f = f'$. 
 
\noindent (We recall that the notion of homotopy employed here means that the two maps are homotopic via equivalences in 
$\C$.)
 \item[(2)] Given a pair of maps $f, f': \Delta^m \to \C$, where $m > n$, if $f_{|\partial \Delta^m} = f'_{|\partial \Delta^m}$, then $f = f'$.
\end{itemize}
These conditions say that $\C$ has no morphisms in degrees $> n$ and any two morphisms in degree $n$ agree if they are equivalent. The conditions can be equivalently expressed as follows: $\C$ is an $n$-category, $n \geq 1$, if for every diagram 
$$
\xymatrix{
\Lambda^m_i \ar[r] \ar[d] & \C \\
\Delta^m \ar@{-->}[ru]
}
$$
where $m > n$ and $0 < i < m$, there exists a \emph{unique} dotted arrow which makes the diagram commutative \cite[Proposition 2.3.4.9]{LurHTT}. Using an inductive argument (see \cite[Proposition 2.3.4.7]{LurHTT}), it can also be shown that the conditions (1) and (2) together are equivalent to:
\begin{itemize}
\item[(3)] Given a simplicial set $K$ and maps $f,f': K \to \C$ such that $f_{| \sk_n(K)}$ and $f'_{|\sk_n(K)}$ are homotopic 
relative to $\sk_{n-1}(K)$, then $f = f'$.
\end{itemize}
An important immediate consequence of (3) is that for every $n$-category $\C$, the $\infty$-category $\Fun(K, \C)$ is again an $n$-category for any simplicial set $K$ \cite[Corollary 2.3.4.8]{LurHTT}. 

\begin{example}
The only $(-1)$-categories up to isomorphism are $\varnothing$ and $\Delta^0$. An $\infty$-category is a $0$-category if and only if it is isomorphic to  
(the nerve of) a poset. $1$-categories are up to isomorphism (nerves of) ordinary categories. See \cite[Examples 2.3.4.2--2.3.4.3, Proposition 2.3.4.5]{LurHTT}. 
\end{example}

The property of being an $n$-category is not invariant under equivalences of $\infty$-categories. The following proposition gives a characterization of the invariant property that an 
$\infty$-category is equivalent to an $n$-category. We recall that an $\infty$-groupoid (= Kan complex) $X$ is \emph{n-truncated}, where $n \geq -1$, if $X$ has vanishing homotopy groups in degrees $> n$. We say that $X$ is \emph{(-2)-truncated} if $X$ is contractible. 

\begin{proposition} \label{characterization}
Let $\C$ be an $\infty$-category and let $n \geq -1$ be an integer. Then $\C$ is equivalent to an $n$-category if and only if $\Map_{\C}(x,y)$ is $(n-1)$-truncated for every pair of objects $x,y \in \C$. 
\end{proposition}
\begin{proof}
See \cite[Proposition 2.3.4.18]{LurHTT}.
\end{proof}

\subsection{Homotopy $n$-categories} Let $\C$ be an $\infty$-category and let $n \geq 1$ be an integer. We recall from \cite{LurHTT} the construction of the homotopy $n$-category $\h_n \C$ of $\C$. Given a simplicial set $K$, we denote by $[K, \C]_n$ the set of maps 
$$\sk_n(K) \to \C$$
which extend to $\sk_{n+1}(K)$. Two elements $f,g \in [K, \C]_n$ are called equivalent, 
denoted $f \sim g$, if the maps $f,g: \sk_n(K) \to \C$ are homotopic relative to $\sk_{n-1}(K)$. The equivalence classes of such maps for $K = \Delta^m$ define the $m$-simplices of a simplicial set $\h_n \C$,  i.e.,
$$(\h_n \C)_m  := [\Delta^m, \C]_n / \sim.$$
Clearly an $m$-simplex of $\C$ defines an $m$-simplex in $\h_n \C$, so we have a canonical map $\gamma_n \colon \C \to \h_n \C$. 
Note that this map is a bijection in simplicial degrees $< n$ and surjective in degrees $n$ and $n+1$.  

\medskip

The following proposition summarises some of the basic properties of this construction.

\begin{proposition} \label{properties_of_htpy_cats} Let $\C$ be an $\infty$-category and $n \geq 1$.
\begin{itemize}
 \item[(a)] The set of maps $K \to \h_n \C$ is in natural bijection with the set $[K, \C ]_n / \sim$. 
 \item[(b)] $\h_n \C$ is an $n$-category. In particular, it is an $\infty$-category. 
 \item[(c)] $\C$ is an $n$-category if and only if the map $\gamma_n \colon \C \to \h_n \C$ is an isomorphism. 
 \item[(d)] Let $\mathscr{D}$ be an $n$-category. Then the restriction functor along $\gamma_n \colon \C \to \h_n \C$,
$$\gamma_n^* \colon \Fun(\h_n \C, \mathscr{D}) \to \Fun(\C, \mathscr{D}),$$
is an isomorphism.
\end{itemize}
\end{proposition}
\begin{proof}
See \cite[Proposition 2.3.4.12]{LurHTT}.
\end{proof}

\begin{example}
For $n = 1$, the $1$-category $\h_{1} \C$ is isomorphic to the (nerve of the) usual homotopy category of $\C$.
\end{example}

\begin{remark} 
For an $\infty$-category $\C$, the homotopy $0$-category $\h_0 \C$ can be described in the following way. For $x, y \in \C$, we write $x \leq y$ if $\Map_{\C}(x, y)$ is non-empty. This defines a reflexive and transitive relation. We say that two objects $x$ and $y$ are equivalent if $x \leq y$ and 
$y \leq x$. Then the relation $\leq$ descends to a partial order on the set of equivalence classes of objects in $\C$. The homotopy $0$-category 
$\h_0 \C$ is isomorphic to the nerve of this poset. We will often ignore the case $n = 0$ and focus on the homotopy $n$-categories for $n \geq 1$. 
\end{remark}

\begin{proposition} \label{cat-fib}
Let $\C$ be an $\infty$-category and $n \geq 1$. The functor $\gamma_n \colon \C \to \h_n \C$ is a categorical fibration between $\infty$-categories. In addition, for every lifting problem
$$
\xymatrix{
\partial \Delta^m \ar[r]^u \ar[d] & \C \ar[d]^{\gamma_n} \\
\Delta^m \ar[r]_{\sigma} \ar@{-->}[ur] & \h_n \C
}
$$
where $m \leq n + 1$ (resp. $m < n$), there is a (unique) filler $\Delta^m \to \C$ which makes the diagram commutative. 
\end{proposition}
\begin{proof}
Clearly for any object $c$ in $\C$ and any equivalence $f \colon c \to c'$ in $\h_n \C$, we may find a lift $\tilde{f} \colon c \to c'$ of $f$ in $\C$ (uniquely if $n>1$), 
which is again an equivalence. Then we need to show that $\gamma_n$ is an inner fibration. Consider a lifting problem
$$
\xymatrix{
\Lambda^m_i \ar[d]_j \ar[r]^u & \C \ar[d]^{\gamma_n} \\
\Delta^m \ar[r]_{\sigma} & \h_n \C
}
$$
where $0 < i < m$. For $m \leq n$, there is a diagonal filler $\Delta^m \to \C$ because $\gamma_n$ is a bijection in simplicial degrees $< n$ and surjective on $n$-simplices. 
For $m > n$, there is a map $v \colon \Delta^m \to \C$ which extends $u$, because $\C$ is an $\infty$-category. We claim that $v$ defines a diagonal filler for the diagram. To see this, note that an extension of $\gamma_n u$ along $j$ is unique up to homotopy relative to $\Lambda^m_i$ ($\supseteq \sk_{n-1} \Delta^m$), since $\h_n \C$ is an $\infty$-category. In particular, the maps $\sigma$ and $\gamma_n v$  are homotopic relative to $\sk_{n-1} \Delta^m$. Then the result follows because the $n$-category $\h_n \C$ satisfies condition (3) (see Subsection \ref{2.1}). Therefore $\gamma_n$ is an inner fibration and this completes the proof of the first claim. 

The second claim for $m \leq n$ follows again from the fact that $\gamma_n$ is bijective in simplicial degrees $< n$ and surjective in degree $n$. 
For $m = n+1$, we may find a map $\sigma' \colon \Delta^{n+1} \to \C$ such that $\gamma_n \sigma' = \sigma$, since $\gamma_n$ is surjective on $(n+1)$-simplices. The maps 
$$u, \sigma'_{| \partial \Delta^{n+1}} \colon \partial \Delta^{n+1} \to \C$$
become equal after postcomposition with $\gamma_n$. By Proposition \ref{properties_of_htpy_cats}(a), this means that they are homotopic (in the sense of the Joyal model category) relative to $\sk_{n-1}(\partial \Delta^{n+1})$. Using that the inclusion $\partial \Delta^{n+1} \subset \Delta^{n+1}$ is a cofibration, this homotopy can be extended to a homotopy on $\Delta^{n+1}$ between $\sigma'$ and a map $v \colon \Delta^{n+1} \to \C$ such that $v_{| \partial \Delta^{n+1}} = u$. This new homotopy is still constant on $\sk_{n-1}(\partial \Delta^{n+1}) = \sk_{n-1}(\Delta^{n+1})$. Therefore, given that $\sigma'$ and $v$ are homotopic relative to $\sk_{n-1}(\Delta^{n+1})$, it follows that $\gamma_n v = \gamma_n \sigma' = \sigma$; this shows that $v$ defines a diagonal filler for the diagram. 
\end{proof}

\begin{corollary} \label{section}
Let $\C$ be an $\infty$-category and let $n \geq 1$ be an integer. There is a (non-canonical) map 
$$\epsilon \colon \sk_{n+1} \h_n \C \to \sk_{n+1} \C$$
such that the following diagram commutes:
$$
\xymatrix{
\sk_{n+1} \C \ar[r] & \C \ar[d]^{\gamma_n} \\
\sk_{n+1} \h_n \C \ar[u]^{\epsilon} \ar[r] & \h_n \C
}
$$  
where the horizontal maps are the canonical inclusions. 
\end{corollary}
\begin{proof}
We have a diagram as follows:
$$
\xymatrix{
\sk_{n-1} \C \ar[d]^{\cong} \ar[r] & \sk_n \C \ar@{->>}[d] \ar[r] & \sk_{n+1} \C \ar[r] \ar@{->>}[d] & \C \ar[d]^{\gamma_n} \\
\sk_{n-1} \h_n \C \ar[r] & \sk_n \h_n \C \ar[r] & \sk_{n+1} \h_n \C \ar[r] & \h_n \C.
}
$$  
We may choose a section $\epsilon' \colon \sk_n \h_n \C \to \sk_n \C$ -- uniquely up to equivalence. We claim that this section can be 
extended further to a section $\epsilon$ as required. Let $\sigma \colon \Delta^{n+1} \to \h_n\C$ denote a nondegenerate $(n+1)$-simplex 
in $\h_n \C$. Then we consider the composite map
$$u \colon \partial \Delta^{n+1} = \sk_n\Delta^{n+1} \xrightarrow{\sk_n(\sigma)} \sk_n \h_n \C \xrightarrow{\epsilon'} \sk_n \C \to \C$$
and the commutative diagram:
$$
\xymatrix{
\partial \Delta^{n+1} \ar[d] \ar[r]^u & \C \ar[d]^{\gamma_n} \\
\Delta^{n+1} \ar[r]_{\sigma} & \h_n \C.
}
$$
By Proposition \ref{cat-fib}, there is a diagonal filler $\tau \colon \Delta^{n+1} \to \C$, that is, an $(n+1)$-simplex of $\C$, which makes 
the diagram commutative. 
We set $\epsilon(\sigma) := \tau$. Repeating this process for each $\sigma$, we obtain the required extension $\epsilon \colon \sk_{n+1} \h_n \C \to \sk_{n+1} \C$. 
\end{proof}

\begin{example}
Let $\C$ be an $\infty$-category. The functor $\gamma_1 \colon \C \to \h_1 \C$ is bijective on objects, so there is a unique section $\sk_0 \h_1 \C \to \sk_0 \C$. By making choices 
of morphisms, one from each homotopy class, this map extends to a section $\sk_1 \h_1 \C \to \sk_1 \C$. The last map extends further to a 
section $\sk_2 \h_1 \C \to \sk_2 \C$ by making (non-canonical) choices of homotopies for compositions. 
\end{example} 

\noindent The functor $\h_n(-)$ preserves categorical equivalences of $\infty$-categories. Using Proposition \ref{properties_of_htpy_cats}, it follows that 
there is a tower of $\infty$-categories:
$$ \C  \to \cdots \to \h_n \C \to \h_{n-1} \C \to \cdots \to \h_{1} \C.$$
By construction, the canonical map
$$\C \longrightarrow \mathrm{lim}(\cdots \to \h_n \C \to \h_{n-1} \C \to \cdots \to \h_1 \C)$$
is an isomorphism, and by Proposition \ref{cat-fib}, this inverse limit defines also a homotopy limit in the Joyal model structure.

\begin{example} \label{truncation_groupoid}
As a consequence of Proposition \ref{characterization}, an $\infty$-groupoid $X$ is categorically equivalent to an $n$-category if and only if it is $n$-truncated. For example, a Kan complex is equivalent to a $0$-category if and only if it is homotopically discrete and to an $1$-category if and only if it is equivalent 
to the nerve of a groupoid. Given an $\infty$-groupoid $X$, the canonical tower of $\infty$-groupoids:
$$X \to \cdots \h_n X \to \h_{n-1} X \to \cdots \to \h_1 X \to \pi_0 X$$   
models the Postnikov tower of $X$ and the map $X \to \h_n X$ is $(n+1)$-connected (i.e., for every $x \in X$, the map $\pi_k(X, x) \to \pi_k(\h_n X, x)$ is a bijection for $k \leq n$ 
and surjective for $k = n+1$.)
\end{example}

\begin{remark} \label{truncation_cat}
An $n$-category $\C$ is $n$-truncated in the $\infty$-category of $\infty$-categories, that is, the $\infty$-groupoid $\Map(K, \C)$ is 
$n$-truncated for any $\infty$-category $K$ (see \cite[5.5.6]{LurHTT} for the definition and properties of truncated objects in an $\infty$-category). 
To see this, recall that $\Fun(K, \C)$ is again an $n$-category and then apply Proposition \ref{characterization}. But 
an $n$-truncated $\infty$-category $\C$ is not equivalent to an $n$-category in general, so the analogue of Example \ref{truncation_groupoid} fails for general $\infty$-categories.  
An $\infty$-category $\C$ is $n$-truncated if and only if $\C$ is equivalent to an $(n+1)$-category and the maximal $\infty$-subgroupoid $\C^{\simeq} \subseteq \C$ is $n$-truncated. Indeed, given an $n$-truncated $\infty$-category $\C$, then $\C^{\simeq} \simeq \Map(\Delta^0, \C)$ is $n$-truncated. Moreover, since $n$-truncated objects are closed under limits, it follows that $\Map_{\C}(x, y)$ is $n$-truncated for 
every $x, y \in \C$, using the fact that there is a pullback in the $\infty$-category of $\infty$-categories:
$$
\xymatrix{
\Map_{\C}(x, y) \ar[d] \ar[r] & \C^{\Delta^1} \ar[d] \\
\Delta^0 \ar[r]^{(x, y)} & \C \times \C.
}
$$
(Using the definition of $n$-truncated objects, we see that $\C^{\Delta^1}$ is again $n$-truncated.) Conversely, if $\C$ is an $(n+1)$-category and $\C^{\simeq}$ is $n$-truncated, then it is possible to show that $\Map(\Delta^k, \C)$ is $n$-truncated by induction on $k \geq 0$, 
from which it follows that $\C$ is $n$-truncated. (I am grateful to Hoang Kim Nguyen for interesting discussions related to this remark.) 
\end{remark}

Let $\catinf$ denote the category of $\infty$-categories, regarded as enriched in $\infty$-categories, and let $\catn$ denote the full subcategory of $\catinf$ which 
is spanned by $n$-categories. 

\begin{proposition}
Let $\C$ and $\D$ be $\infty$-categories and let $n \geq 1$ be an integer. 
\begin{itemize}
\item[(a)] The natural map $\h_n (\C \times \D) \xrightarrow{\cong} \h_n \C \times \h_n \D$ is an isomorphism. 
\item[(b)] There is a functor 
$$\h_{n}^{\C, \D} \colon \Fun(\C, \D) \to \Fun(\h_n \C, \h_n \D)$$
which is natural in $\C$ and $\D$. In particular, $\h_n \colon \catinf \to \catn$ is an enriched functor. 
\end{itemize}
\end{proposition}
\begin{proof}
(a) follows directly from the definition of $\h_n$. For (b), we define the functor $\h_n^{\C, \D}$ as follows: a $k$-simplex $F \colon \C \times \Delta^k \to \D$ is sent to the composite
$$\h_n \C \times \Delta^k \cong \h_n \C \times \h_n \Delta^k \cong \h_n (\C \times \Delta^k) \xrightarrow{\h_n F} \h_n \D.$$
The functor $\h_n^{\C, \D}$ is natural in $\C$ and $\D$ and turns $\h_n$ into an enriched functor. 
\end{proof}

\section{Higher weak colimits} \label{sec:higher_weak_colimits}

\subsection{Basic definitions and properties}  It is well known that homotopy categories do not admit small (co)limits in general, even when the underlying $\infty$-category has small (co)limits. On the other hand, if the $\infty$-category $\C$ admits pushouts (resp. coproducts), for example,  then the homotopy category $\h_1 \C$ admits weak pushouts (resp. coproducts), which are induced from pushouts (resp. coproducts) in $\C$. Moreover, if $\C$ admits small colimits, then $\h_1 \C$ admits small weak colimits -- which may or may not be induced from $\C$. 
These observations suggest the following questions: does $\h_n \C$, $n > 1$, have in some sense more or \emph{better} (co)limits than the homotopy category, and how do these compare with (co)limits in $\C$?

\medskip 

In this section we introduce and study a notion of higher weak (co)limit in the context of $\infty$-categories. This is both a higher categorical version
of the classical notion of weak (co)limit and a weak version of the higher categorical notion of (co)limit. We will mainly focus on the properties of higher weak (co)limits in the context of higher homotopy categories. We also refer to \cite{more-gaft} for further properties and applications of higher weak (co)limits.  

\medskip

We will restrict to higher weak colimits since the corresponding definitions and results about higher weak limits are obtained dually.  First we recall that an $\infty$-groupoid (= Kan complex) $X$ is $k$-\emph{connected}, for some $k \geq -1$, if it is non-empty and $\pi_i(X, x) \cong 0$ for every $x \in X$ and $i \leq k$. For example, $X$ is $(-1)$-connected (resp. $0$-connected, $\infty$-connected) if it is non-empty (resp. connected, contractible).  

\medskip

We begin with the definition of a weakly initial object. Fix $t \in \mathbb{Z}_{\geq 0} \cup \{\infty\}$.  

\begin{definition}
An object $x$ of an $\infty$-category $\C$ is \emph{weakly initial of order} $t$ if the mapping space $\Map_{\C}(x, y)$ is $(t-1)$-connected for every object $y \in \C$.
\end{definition}

\begin{example} \label{weakly-initial-example}
If $\C$ is an ordinary category, a weakly initial object $x \in \C$ of order $0$ is a weakly initial object in the classical sense. For a general $n$-category 
$\C$, a weakly initial object of order $n$ is an initial object.
\end{example}

\begin{proposition}
Let $\C$ be an $\infty$-category and $x \in \C$. The following are equivalent: 
\begin{itemize}
\item[(1)] $x$ is weakly initial in $\C$ of order $t$.
\item[(2)] $x$ is weakly initial in $\h_n \C$ of order $t$ for any $n > t$.
\item[(3)] $x$ is initial in $\h_{t} \C$.
\end{itemize}
These imply:
\begin{itemize}
\item[(4)] $x$ is initial in $\h_n \C$ for any $n < t$. 
\end{itemize}
\end{proposition}
\begin{proof}
This follows from the fact that the functor $\gamma_n \colon \C \to \h_n \C$ restricts to an $n$-connected map $\Map_{\C}(x, y) \to \Map_{\h_n \C}(x, y)$ 
for every $x, y \in \C$.
\end{proof}

\begin{proposition} \label{Kan_complex_weakly_initial}
Let $\C$ be an $\infty$-category and let $t > 0$. The full subcategory $\C'$ of $\C$ which is spanned by the weakly initial objects of order $t$ is either empty or 
a $t$-connected $\infty$-groupoid. 
\end{proposition}
\begin{proof}
Suppose that the full subcategory $\C'$ is non-empty. Then the mapping spaces of $\C'$ are $(t-1)$-connected, where $t-1 \geq 0$. It follows that every morphism in $\C'$ is an equivalence, therefore $\C'$ is an $\infty$-groupoid. 
\end{proof}

\begin{remark}
The full subcategory $\C'$ of weakly initial objects in $\C$ of order 0 is not an $\infty$-groupoid in general. In this case, we only have that $\Map_{\C}(x, y)$ 
is non-empty for every $x, y \in \C'$.
\end{remark} 

\begin{definition}
Let $\C$ be an $\infty$-category, $K$ a simplicial set, and let $F \colon K \to \C$ be a $K$-diagram in $\C$. 
A weakly initial object $G \in \C_{F/}$ of order $t$ is called a \emph{weak colimit of $F$ of order} $t$.   
\end{definition}

\begin{example}
If $\C$ is an $n$-category and $G \colon K^{\triangleright} \to \C$ is a weak colimit of $F = G_{|K} \colon K \to \C$ order $t \geq n$, then $G$ is a colimit diagram. This follows from Example \ref{weakly-initial-example} using the fact that $\C_{F/}$ is again an $n$-category (see \cite[Corollary 2.3.4.10]{LurHTT}). In particular, a weak colimit of order $\infty$ is a colimit diagram.  If $\C$ is an ordinary category, then a weak colimit of order $0$ is a weak colimit diagram in the classical sense. 
\end{example}

\begin{remark} \label{difference-of-higher-weak-colimits}
There is an important difference between weak colimits of order $0$ and weak colimits of order $> 0$: as a consequence of Proposition \ref{Kan_complex_weakly_initial}, any two weak colimits 
of $F$ of order $> 0$ are equivalent. In particular, if $G \in \C_{F/}$ is a weak colimit of order $t > 0$ and $G' \in \C_{F/}$ is a weak colimit of order $>0$, then $G'$ is also a weak colimit of 
order $t$. 
\end{remark}

The following proposition gives an alternative characterization of higher weak colimits following the analogous characterization for colimits in \cite[Lemma 4.2.4.3]{LurHTT}. 

\begin{proposition} \label{alternative_def}
Let $\C$ be an $\infty$-category, $K$ a simplicial set, and let $G \colon K^{\triangleright} \to \C$ be a diagram with cone object $x \in \C$ . Then $G$ is a weak colimit 
of $F = G_{| K}$ of order $t$ if and only if the canonical restriction map: 
$$\Map_{\C}(x, y) \simeq \Map_{\C^{K^{\triangleright}}}(G, c_y) \to \Map_{\C^K}(F, c_y)$$
is $t$-connected for every $y \in \C$, where $c_y$ denotes respectively the constant diagram at $y \in \C$. 
\end{proposition}
\begin{proof}
The fiber of the restriction map over $F' \colon K^{\triangleright} \to \C$ with cone object $y$ is identified with $\Map_{\C_{F/}}(G, F')$ (see the proof of \cite[Lemma 4.2.4.3]{LurHTT}).
\end{proof}

The basic rules for the manipulation of higher weak colimits can be established similarly as for colimits. The following procedure shows that higher weak colimits can be computed iteratively, exactly like colimits, but with the difference that the order of the weak colimit may decrease with each iteration. 

\begin{proposition} \label{iterated_weak_colimits}
Let $\C$ be an $\infty$-category and let $K = K_1 \cup_{K_0} K_2$ be a simplicial set where $K_0 \subseteq K_1$ is a simplicial subset. Let $F \colon K \to \C$ be a diagram and denote its restrictions by $F_i := F_{|K_i }$, for $i = 0,1,2$. Suppose that $G_i \colon K_i^{\triangleright} \to \C$ is a weak colimit  of $F_i$ of order $t_i \geq 0$, for $i = 0,1,2$.  
\begin{itemize}
\item[(a)] There are morphisms $G_0 \to G_1{}_{|K_0^{\triangleright}}$ and $G_0 \to G_2{}_{|K_0^{\triangleright}}$ in $\C_{F_0/}$. These together with $G_1$ and $G_2$ 
determine a diagram in $\C$ as follows, 
$$H \colon K_1^{\triangleright} \cup_{K_0  \ast \Delta^{\{1\}}} (K_0 \ast \Delta^1) \cup_{K_0 \ast \Delta^{\{0\}}} (K_0 \ast \Delta^1) 
\cup_{K_0 \ast \Delta^{\{1\}}}  K_2^{\triangleright}  \to \C.$$
\item[(b)] Let $H_{\ulcorner} \colon \Delta^1 \cup_{\Delta^0} \Delta^1 \to \C$ be the restriction of $H$ to the cone objects. Suppose that 
$H' \colon \Delta^1 \times \Delta^1 \to \C$ is a weak pushout of $H_{\ulcorner}$ of order $k$. Then $H'$ determines a cocone $G \colon K^{\triangleright} \to \C$ over $F$, with the 
same cone object as $H'$, which is a weak colimit of $F$ of order $\ell : = \mathrm{min}(k, t_1, t_0 -1, t_2)$.
\end{itemize} 
\end{proposition}
\begin{proof}
(a) follows directly from the properties of higher weak colimits. For (b), we first explain the construction of the cocone $G \colon K^{\triangleright} \to \C$. The functor $H'$ is represented by a diagram
$$
\xymatrix{
x_0 \ar[d]_{v} \ar[r]^u & x_1 \ar[d]^f \\
x_2 \ar[r]_g & y
}
$$
where $x_i$ is the cone object of $G_i$, and the morphisms $u$ and $v$ are given respectively by the morphisms $G_0 \to G_1{}_{|K_0^{\triangleright}}$ and $G_0 \to G_2{}_{|K_0^{\triangleright}}$ in $\C_{F_0/}$. The morphisms $f$ and $g$ produce two new 
cocones (essentially uniquely): $G'_1 \colon K_1^{\triangleright} \to \C$ over $F_1$, and $G'_2 \colon K_2^{\triangleright} \to \C$ over $F_2$, with common cone object $y$. The restrictions $G'_1{}_{|K_0^{\triangleright}}$ and 
$G'_2{}_{|K_0^{\triangleright}}$ are equivalent as cocones over $F_0$. We may then extend $G'_2{}_{|K_0^{\triangleright}}$ in an essentially unique way
to a new cocone $G''_1 \colon K_1^{\triangleright} \to \C$ over $F_1$, which is equivalent to $G'_1$. The resulting cocones 
$$G'_2{}_{|K_0^{\triangleright}} \colon K_0^{\triangleright} \to \C, \ G''_1\colon K_1^{\triangleright} \to \C,  \text{and} \ G'_2 \colon K_2^{\triangleright} \to \C$$ assemble to  define the required cocone $G \colon K^{\triangleright} \to \C$. Then the claim in (b) is shown by applying Proposition \ref{alternative_def}, first for weak pushouts and then for $K_i$-diagrams, and using the fact that $\C^K$ is the homotopy pullback (in the Joyal model category) of the diagram of $\infty$-categories ($\C^{K_1} \to \C^{K_0} \leftarrow \C^{K_2}$), so its mapping spaces can also be identified with the corresponding (homotopy) pullbacks of mapping spaces.
\end{proof}

\begin{example}
Let $\C$ be (the nerve of) an ordinary category that admits small coproducts and weak pushouts. By Proposition \ref{iterated_weak_colimits}, every diagram $F \colon K \to \C$, where $K$ is 1-dimensional, admits a weak colimit (of order 0). Now suppose that $F \colon I \to \C$ is a diagram where $I$ is (the nerve of) an arbitrary ordinary small category. Since $\C$ is an $1$-category, a cocone $G \colon I^{\triangleright} \to \C$ over $F$ is determined uniquely by its restriction to a 
cocone $G' \colon (\sk_1 I)^{\triangleright} \to \C$ over $F_{|\sk_1 I}$, and similarly for morphisms between cocones. As a consequence, we may deduce the well known fact that $\C$ has weak $I$-colimits.  
\end{example}

\begin{example}
Let $T \subset \Delta^1 \times \Delta^2$ be the full subcategory spanned by the objects $(0, i)$, for $i = 0,1,2$, and $(1,0)$. Let $\C$ be an $\infty$-category with weak pushouts of 
order $t$ and let $F \colon T \to \C$ be a $T$-diagram in $\C$. Write $T = T_1 \cup_{T_0} T_2$ where $T_1$ is spanned by $(0,i)$, $i = 0,1$, and $(1,0)$, $T_2$ is spanned by $(0, i)$, 
$i = 1,2$, and $T_0 = \{(0,1)\}$. Using Proposition \ref{iterated_weak_colimits}, we may compute a weak colimit of $F$ of order $t$ in terms of iterated weak pushouts of order $t$. 
\end{example}

\subsection{Homotopy categories and (co)limits}
Let $\C$ be an $\infty$-category and let $K$ be a simplicial set. By the universal property of $\h_n(-)$, the functor 
$$\Fun(K, \C) \to \Fun(K, \h_n \C),$$
which is given by composition with $\gamma_n \colon \C \to \h_n \C$, factors canonically through the homotopy $n$-category: 
\begin{equation}
\Phi^K_n: \h_n \Fun(K, \C) \to \Fun(K, \h_n \C).
\end{equation}
The comparison between $K$-colimits in $\C$ and in $\h_n \C$ is essentially a question about the properties of the functor $\Phi^K_n$. 
Note that for $n=1$, $\Phi^K_1$ is simply the canonical functor of ordinary categories: $\mathrm{h}_1(\C^K) \to \mathrm{h}_1(\C)^K$. 

\begin{lemma} \label{lem:simplices}
Let $\C$ be an $\infty$-category, $K$ a finite dimensional simplicial set of dimension $d > 0$, and let $n \geq 1$ be an integer. The functor 
$$\Phi^K_n : \h_n \Fun(K, \C) \to \Fun(K, \h_n \C)$$ 
satisfies the following:
\begin{itemize}
\item[(a)] $\Phi^K_n$ is a bijection in simplicial degrees $< n - d$.
\item[(b)] $\Phi^K_n$ is surjective in simplicial degree $n-d$; it identifies $(n-d)$-simplices which are homotopic relative to the $(n-1)$-skeleton of $\Delta^{n-d} \times K$.
\item[(c)] $\Phi^K_n$ is surjective in simplicial degree $n-d+1$. 
\end{itemize}
\end{lemma}
\begin{proof}
The $m$-simplices of $\Fun(K, \h_n \C)$ are equivalence classes of maps
$$\sk_n(\Delta^m \times K) \to \C$$
that extend to $\sk_{n+1}(\Delta^m \times K)$. On the other hand, the $m$-simplices of $\h_n \Fun(K, \C)$ are equivalence classes of maps 
$$\sk_n(\Delta^m) \times K \to \C$$
that extend to $\sk_{n+1}(\Delta^m) \times K$. The functor $\Phi^K_n$ is induced by the canonical map 
$$\sk_n(\Delta^m \times K) \to \sk_n(\Delta^m) \times K$$
which is an isomorphism if $d \leq \mathrm{max}(n - m, 0)$. This shows that the map $\Phi^K_n$ is 
surjective in simplicial degrees $\leq n-d + 1$. Similarly, the map 
$$\sk_{n-1}(\Delta^m \times K) \to \sk_{n-1}(\Delta^m) \times K$$
is an isomorphism if $d \leq \mathrm{max}(n - m -1, 0)$, so the two equivalence relations agree for $m < n -d$. 
\end{proof}

\begin{remark} \label{special case}
The case $d = 0$ is both special and essentially trivial, since the functor $\Phi^K_n$ is an isomorphism in this case.  
\end{remark}

\begin{proposition} \label{lifting-prop-Phi}
Let $\C$ be an $\infty$-category, $K$ a finite dimensional simplicial set of dimension $d > 0$, and let $n \geq 1$ be an integer. Then for every lifting problem
$$
\xymatrix{
\partial \Delta^m \ar[r]^u \ar[d] & \h_n(\C^K) \ar[d]^{\Phi^K_n} \\
\Delta^m \ar[r]_{\sigma} \ar@{-->}[ur] & (\h_n \C)^K
}
$$
where $m \leq n-d+1$ (resp. $m < n - d$), there is a (unique) filler $\Delta^m \to \h_n (\C^K)$ which makes the diagram commutative. 
\end{proposition}
\begin{proof}
The case $m < n-d + 1$ is a direct consequence of Lemma \ref{lem:simplices}. For $m = n -d +1 \leq n$, Lemma \ref{lem:simplices} shows that there is a lift $\tau \colon \Delta^m \to \h_n(\C^K)$ of $\sigma$, represented by a map $\tau' \colon \Delta^m \to \C^K$. The maps 
$$u,  \gamma_n \circ \tau'_{|\partial \Delta^m} \colon \partial \Delta^m \to \h_n(\C^K)$$ become equal after composition with $\Phi^K_n$. Moreover,  note that $u$ corresponds uniquely to a map $u \colon \partial \Delta^m \to \C^K$. The maps $\Phi^K_n \circ u$ and $\Phi^K_n \circ (\gamma_n \tau'_{|\partial \Delta^m})$ correspond to the equivalence classes of maps (using here the same notation for the adjoint maps)
$$u, \tau'_{| \partial \Delta^m \times K} \colon \partial \Delta^m \times K \to \C,$$ 
so these last two maps are homotopic relative to $\sk_{n-1}(\partial \Delta^m \times K)$ (see Proposition \ref{properties_of_htpy_cats}). Let $$H \colon (\partial \Delta^m \times K) \times J \to \C$$ be a homotopy from $\tau'_{| \partial \Delta^m \times K}$ to $u$ relative to the subspace $\sk_{n-1}(\partial \Delta^m \times K)$, where $J = N(0 \rightleftarrows 1)$ denotes the Joyal interval object. Using the fact that 
$$\partial \Delta^m  \times K \cup_{\sk_{n-1}(\partial \Delta^m \times K)} \sk_{n-1}(\Delta^m \times K) \subset \Delta^m \times K$$ is a cofibration, 
this homotopy can be extended to a homotopy 
$$H' \colon (\Delta^m \times K) \times J \to \C$$
where $\tau' = H'_{|\Delta^m\times K \times \{0\}}$ and $H'$ is constant on $\sk_{n-1}(\Delta^m \times K)$. Then the map 
$$\widetilde{\tau} \colon = H'_{| \Delta^m \times K \times \{1\}} \colon \Delta^m \times K \to \C$$ 
extends $u$. Moreover, the (adjoint) map $\widetilde{\tau} \colon \Delta^m \to \C^K$ defines an element in $\h_n(\C)^K_m$, after composition with $\gamma_n$, which agrees with 
$\Phi^K_n \circ \tau = \sigma$, so $\widetilde{\tau}$ represents a diagonal filler, as required.  
\end{proof}

As a consequence of Proposition \ref{lifting-prop-Phi}, we obtain the following result about the comparison between the $n$-categories $\h_n(\C^K)$ and $\h_n(\C)^K$. 

\begin{corollary} \label{n-d_homotopy_cat}
Let $\C$ be an $\infty$-category, $K$ a finite dimensional simplicial set of dimension $d > 0$, and let $n \geq d$ be an integer. The functor $\Phi^K_n \colon \h_n(\C^K) \to 
(\h_n \C)^K$ is essentially surjective and for every pair of objects $F, G$ in $\h_n(\C^K)$, the induced map between mapping spaces
$$\Map_{\h_n(\C^K)}(F, G) \longrightarrow \Map_{(\h_n \C)^K}(\Phi^K_n(F), \Phi^K_n(G))$$
is $(n-d)$-connected. As a consequence, $\Phi^K_n$ induces an equivalence:
\begin{equation} \label{n-d_homotopy_cat_eqn} 
\h_{n-d}\big(\Fun(K,\C)\big) \simeq \h_{n-d}\big(\Fun(K, \h_n \C)\big). 
\end{equation}
\end{corollary}

\medskip

Now assume that $\C$ is an $\infty$-category which admits $K$-colimits, where $K$ is a simplicial set. We have colimit-functors:
$$
\xymatrix{
\C^K \ar[rr]^{\mathrm{colim}_K} \ar[d]_{\gamma_n} && \C \ar[d]^{\gamma_n} \\
\h_n( \C^K) \ar[rr]^(.5){\h_n(\mathrm{colim}_K)} && \h_n \C.
}
$$
Assuming also that $\C$ and $K$ are as in Corollary \ref{n-d_homotopy_cat}, and passing to the homotopy $(n-d)$-categories as 
in \eqref{n-d_homotopy_cat_eqn}, we obtain the following corollary.

\begin{corollary} \label{n-d_homotopy_cat2}
Let $\C$ be an $\infty$-category, $K$ a finite dimensional simplicial set of dimension $d > 0$, and let $n \geq d$ be an integer. Suppose that $\C$ admits $K$-colimits. Then there is an adjoint pair
$$\h_{n-d}(\colim_K) \colon \h_{n-d}\big(\Fun(K, \h_n \C)\big) \rightleftarrows \h_{n-d} (\C) \colon \h_{n-d}(c)$$
where $c$ denotes the constant $K$-diagram functor.
\end{corollary}

\begin{proposition} \label{preservation of weak colimits}
Let $\C$ be an $\infty$-category with weak $K$-colimits of order $k$, where $K$ is a simplicial set of dimension $d > 0$, and let $n \geq 1$ be an integer. Then $\h_n \C$ has weak $K$-colimits of order $\ell = \mathrm{min}(n-d, k)$. Moreover, the functor $\gamma_n \colon \C \to \h_n \C$ sends weak $K$-colimits of order $k$ to weak $K$-colimits of order $\ell$. 
\end{proposition}
\begin{proof}
We may assume that $n \geq d$ and therefore the functor $\C^K \to (\h_n \C)^K$ is surjective on objects. Then it suffices to prove the second claim. 
Let $G \colon K^{\triangleright} \to \C$ be a weak colimit of $F = G_{|K}$ of order $k$ with cone object $x \in \C$. By Proposition \ref{alternative_def}, it suffices to prove that the canonical map 
$$\Map_{(\h_n\C)^{K^{\triangleright}}}(G, c_y) \to \Map_{(\h_n\C)^K}(F, c_y)$$
is $\ell$-connected for every $y \in \C$.  Note that there is a $k$-connected map: 
$$\Map_{(\h_n \C)^{K^{\triangleright}}}(G, c_y) \simeq \Map_{\h_n \C}(x, y) \simeq \Map_{\h_n(\C^{K^{\triangleright}})}(G, c_y) \to \Map_{\h_n(\C^K)}(F, c_y).$$  
Hence it suffices to show that the canonical map 
$$\Map_{\h_n(\C^K)}(F, c_y) \to \Map_{(\h_n\C)^K}(F, c_y)$$
is $(n-d)$-connected. This follows from Corollary \ref{n-d_homotopy_cat}. (Alternatively, note that the last map is identified with the canonical map from the $(n-1)$-truncation of a $K^{\op}$-limit of $\infty$-groupoids to the $K^{\op}$-limit of the $(n-1)$-truncations of the $\infty$-groupoids:
$$\h_{n-1} \big(\mathrm{lim}_{K^{\op}} \Map_{\C}(F(-), y)\big) \to \mathrm{lim}_{K^{\op}} \h_{n-1} \big(\Map_{\C}(F(-), y)\big).$$ 
An inductive argument on $d$ shows that the map is $(n-d)$-connected.)
\end{proof}

\begin{remark}
A different proof of Proposition \ref{preservation of weak colimits} is also possible using elementary lifting arguments based on Proposition \ref{lifting-prop-Phi}. 
\end{remark}

\begin{corollary} \label{cor_weak_pushouts}
Let $\C$ be an $\infty$-category which admits finite (resp. small) colimits. 
\begin{itemize}
\item[(a)] The homotopy $n$-category $\h_n \C$ admits finite (resp. small) coproducts and weak pushouts of order $n-1$. Moreover, the functor $\gamma_n \colon \C \to \h_n \C$ preserves finite (resp. small) coproducts and sends pushouts in $\C$ to weak pushouts of order $n-1$. 
\item[(b)] Suppose that $\gamma_n \colon \C \to \h_n \C$ preserves finite colimits. Then $\C$ is equivalent to an $n$-category. 
\end{itemize}
\end{corollary}
\begin{proof}
(a) $\h_n \C$ admits finite coproducts by Remark \ref{special case}. The existence and preservation of higher weak pushouts is a consequence of Proposition \ref{preservation of 
weak colimits}. (b) is a consequence of \cite[Corollary 3.3.5]{gaft}. 
\end{proof}

These properties of homotopy $n$-categories of (finitely) cocomplete $\infty$-categories suggest to consider the following class of $n$-categories as a convenient general context for the study of these $n$-categories. 

\begin{definition} \label{weakly-cocomplete-n-cat}
Let $n \geq 1$ be an integer or $n=\infty$. A \emph{(finitely) weakly cocomplete $n$-category} is an $n$-category $\C$ which  admits small (finite) coproducts and weak pushouts of order $n-1$. 
\end{definition}

These $n$-categories are studied further in \cite{more-gaft} in connection with adjoint functor theorems and Brown representability theorems for $n$-categories (extending the results of \cite{gaft}).  Also, in Section \ref{sec:K-theory_higher_homotopy_cats} (which can be read independently of Sections \ref{sec:higher_derivators} and \ref{sec:K-theory_higher_derivators}), we will return to the properties of weak (co)limits in higher homotopy categories, especially, for pointed and stable $\infty$-categories, and use higher weak colimits in order to define $K$-theory for higher homotopy categories. 
 
 \begin{remark} \label{partial-colimits} Corollary \ref{n-d_homotopy_cat2} produces a \emph{truncated $K$-colimit functor} for $\h_n \C$:
$$\h_{n-d}(\colim_K) \colon \h_{n-d}\big(\Fun(K, \h_n \C)\big) \to \h_{n-d}(\C).$$
According to Proposition \ref{lifting-prop-Phi} (cf. Proposition \ref{cat-fib} and Corollary \ref{section}), there is a (non-canonical) section
$$\epsilon \colon \sk_{n-d+1}\big(\Fun(K, \h_n\C)\big) \to \sk_{n-d+1}\big(\h_n(\C^K)\big).$$
(Note that the target is closely related to $\sk_{n-d+1}(\C^K).$) Then we may define the \emph{partial $K$-colimit functor} for $\h_n (\C)$ (which depends on $\C$ and the section $\epsilon$):
$$\sk_{n-d+1} \big(\Fun(K, \h_n \C) \big) \to \h_n(\C)$$
as the composition
$$\sk_{n-d+1}\big(\Fun(K, \h_n\C)\big) \xrightarrow{\epsilon} \sk_{n-d+1}\big(\h_n(\C^K)\big) \to \h_n(\C^K) \xrightarrow{\h_n(\colim_K)} \h_n(\C).$$
\end{remark}
 
\begin{example} 
Let $\C$ be an $\infty$-category which has pushouts and let $K$ denote the (nerve of the) ``corner'' category $\ulcorner : = \Delta^1 \cup_{\Delta^0} \Delta^1$. The functor
$$\Phi^K_1 \colon \h_{1} (\C^K) \to (\h_{1} \C)^K$$
is surjective on objects and full. By Corollary \ref{n-d_homotopy_cat}, the pushout-functor on $\C$ induces a truncated pushout-functor:
$$\h_0(\colim_K) \colon \h_0 \big(\Fun(K, \h_1 \C)\big) \to \h_0(\C)$$
which is left adjoint to the constant diagram functor. Furthermore, we have a map as follows,
\begin{equation} \label{partial_colimit_example}
\sk_0\big(\Fun(K, \h_{ 1} \C)\big)  \to \h_{ 1}(\C)
\end{equation}
which sends $F \colon K \to \h_1 \C$ to the pushout of a choice of a lift $\widetilde{F} \colon K \to \C$. This is simply regarded as a map from 
the \emph{set} of $0$-simplices. Moreover, \eqref{partial_colimit_example} extends further to the $1$-skeleton, but this involves non-canonical 
choices which are not unique even up to homotopy. As explained in Remark \ref{partial-colimits}, an extension of this type can be obtained from 
a section $\epsilon \colon \sk_1 \big(\Fun(K, \h_1\C)\big) \to \sk_1\big(\h_1(\C^K)\big).$
The fact that this process cannot be continued to higher dimensional skeleta is related to the non-functoriality of weak pushouts in the homotopy category.  

More generally, for $n \geq 1$, there is a left adjoint truncated pushout-functor:
$$\h_{n-1}\big(\Fun(K, \h_n \C)\big) \to \h_{n-1}(\C)$$
and partial pushout-functors:
$$\sk_{n}\big(\Fun(K, \h_n \C)\big) \to \h_n(\C)$$
which define weak pushouts in $\h_n(\C)$ of order $n-1$. 
\end{example}

\section{Higher derivators} \label{sec:higher_derivators}

\subsection{Basic definitions and properties} We recall that $\catinf$ denotes the category of $\infty$-categories, regarded as enriched in $\infty$-categories. Let 
$\dia$ denote a full subcategory of $\catinf$ which has the following properties:
\begin{itemize} 
 \item[(Dia 0)] $\dia$ contains the (nerves of) finite posets. 
 \item[(Dia 1)] $\dia$ is closed under finite coproducts and under pullbacks along an inner fibration.
 \item[(Dia 2)] For every $\x \in \dia$ and $x \in \x$, the $\infty$-category $\x_{/x}$ is in $\dia$. 
 \item[(Dia 3)] $\dia$ is closed under passing to the opposite $\infty$-category. 
\end{itemize}
The main examples of such subcategories of $\catinf$ are the following:
\begin{itemize}
\item[(i)] The full subcategory of (nerves of) finite posets. 
\item[(ii)] The full subcategory of (ordinary) finite direct categores $\dirf$. (We recall that an ordinary category $\c$ is called \emph{finite direct} if its nerve is a finite simplicial set.)
\item[(iii)] The full subcategory $\catn \subset \catinf$ of $n$-categories for any $n \geq 1$. 
\item[(iv)] $\catinf$.
\end{itemize}  
We denote by $\dia^{\op}$ the opposite category taken 1-categorically, that is, the enrichment of $\dia^{\op}$ is given by:
$$\underline{\mathrm{Hom}}_{\dia^{\op}}(\x, \y) = \underline{\mathrm{Hom}}_{\dia}(\y, \x) = \Fun(\y, \x).$$

\begin{definition}
An $\infty$-\emph{prederivator} with domain $\dia$ is an enriched functor 
$$\der \colon \dia^{\op} \to \catinf.$$ 
An $\infty$-prederivator $\der$ with domain $\dia$ is an $n$-\emph{prederivator} 
if it factors through the inclusion $\catn \subset \catinf$, that is, $\der$ is an enriched functor 
$$\der \colon \dia^{\op} \to \catn.$$ 
\end{definition}

A \emph{strict} morphism of $\infty$-prederivators is a natural transformation $F \colon \der \to \der'$ between enriched functors. Thus, 
we obtain a category of $\infty$-prederivators, denoted by $\PreDer_{\infty}$, which is enriched in $\infty$-categories. For any $n \geq 1$, there is a full 
subcategory $\PreDer_n \subset \PreDer_{\infty}$ spanned by the $n$-prederivators. In the same way that the classical theory of (pre)derivators is founded 
on a ((2,2)=)2-categorical context, the general theory of $\infty$-prederivators involves an $(\infty, 2)$-categorical context. We point out that it would be
natural to consider also non-strict morphisms (aka. pseudonatural transformations) between $\infty$-prederivators, but these will not be needed in this paper. 

\medskip 

\noindent \textbf{Notation.} Let $\der$ be an $\infty$-prederivator with domain $\dia$ and let $u \colon \x \to \y$ be a functor in $\dia$. We will often denote 
the induced functor $\der(u) \colon \der(\y) \to \der(\x)$ by $u^*$. Moreover, if $i_{\y, y} \colon \Delta^0 \to \y$ is the inclusion of the object $y \in \y$ and $F \in \der(\y)$, we will often denote the object $\der(i_{\y, y})(F)$ in $\der(\Delta^0)$ by $F_{y}$. 

\medskip

\begin{example}
An 1-prederivator with domain $\mathsf{Cat_1}$ is a prederivator in the usual sense \cite{Malt, Groth}. 
\end{example}

\begin{example} \label{represented_der}
Let $\C$ be an $\infty$-category. There is an associated $\infty$-prederivator (with domain $\dia$) defined by $\der^{(\infty)}_\C \colon \dia^{\op} \to \catinf,  \ \x \mapsto \Fun(\x, \C)$.
Moreover, for any $n \geq 1$, 
$$\der^{(n)}_\C \colon \dia^{\op} \to \catn, \ \ \x \mapsto \h_n \big(\Fun(\x, \C)\big),$$
defines an $n$-prederivator. 
\end{example}

\begin{definition} \label{deriv_def}
An $\infty$-prederivator $\der \colon \dia^{\op} \to \catinf$ is a \emph{right} $\infty$-\emph{derivator} if it satisfies the following properties:

\smallskip 

\renewcommand{\theenumi}{(Der \arabic{enumi})}
\renewcommand{\labelenumi}{\theenumi}
\begin{enumerate}
\setcounter{enumi}{0}
\setlength{\itemsep}{.2cm}
\item For every pair of $\infty$-categories $\x$ and $\y$ in $\dia$, the functor induced by the inclusions of the factors to the
coproduct $\x \sqcup \y$,
$$\der(\x \sqcup \y)\rightarrow \der(\x) \times \der(\y),$$
is an equivalence. Moreover, $\der(\varnothing)$ is equivalent to the final $\infty$-category $\Delta^0$.

\item For every $\infty$-category $\x$ in $\dia$, the functor
$$(i_{\x, x}^{*} = \der(i_{\x, x}))_{x \in \x}  \colon \der(\x) \rightarrow \prod_{x \in \x} \der(\Delta^0)$$
is conservative, i.e., it detects equivalences. We recall that $i_{\x, x} \colon \Delta^0 \to \x$ is the functor that corresponds to the object $x \in \x$.

\item For every morphism $u \colon \x \rightarrow \y$ in $\dia$, the functor $u^{*} = \der(u) \colon\der(\y)\rightarrow \der(\x)$ is a right adjoint. We denote a left adjoint of $u^*$ by
$$u_{!}\colon \der(\x)\rightarrow\der(\y).$$

\item Given $u \colon \x \rightarrow \y$ in $\dia$ and $y \in \y$, consider the following pullback diagram in $\dia$,
$$
\xymatrix{
u_{/y} \ar[d]_{p_{u / y}}\ar[r]^{j_{u / y}} & \x \ar[d]^-{u}\\
\y_{/y} \ar[r]_-{q_{\y / y}}& \y. 
}
$$
Then the canonical base change natural transformation:
$$c_{u,y}\colon (p_{u/y})_{!} j^{*}_{u/y} \longrightarrow q^{*}_{\y /y }u_{!}$$
is a natural equivalence of functors.
\end{enumerate} 
\end{definition}

\noindent We define left $\infty$-derivators dually.

\begin{definition}
An $\infty$-prederivator $\der \colon \dia^{\op} \to \catinf$ is a \emph{left} $\infty$-\emph{derivator} if it satisfies (Der1)--(Der2) as stated above together with the following dual versions of (Der3)--(Der4):

\smallskip

\begin{enumerate}
\setlength{\itemsep}{.2cm}
\item[(Der 3)*] For every morphism $u \colon \x \rightarrow \y$ in $\dia$, the functor $u^{*} = \der(u) \colon\der(\y)\rightarrow \der(\x)$ is a left adjoint. 
We denote a right adjoint of $u^*$ by $$u_*\colon \der(\x) \to \der(\y).$$
\item[(Der 4)*] Given $u \colon \x \rightarrow \y$ in $\dia$ and $y \in \y$, consider the following pullback diagram in $\dia$,
$$
\xymatrix{
u_{y/} \ar[d]_{p_{y / u}}\ar[r]^{j_{y / u}} & \x \ar[d]^-{u}\\
\y_{y/} \ar[r]_-{q_{y / \y}}& \y. 
}
$$
Then the canonical base change natural transformation:
$$c'_{u,y} \colon q^{*}_{y / \y }u_{*} \longrightarrow (p_{y/u})_{*} j^{*}_{y/u}$$
is a natural equivalence of functors.  
\end{enumerate}
\end{definition}

\begin{example}
Let $\der \colon \dia^{\op} \to \catinf$ be a left $\infty$-derivator. Then the $\infty$-prederivator
$$\der(-^{\op})^{\op}\colon \dia^{\op} \to \catinf, \ \x \mapsto \der(\x^{\op})^{\op},$$ 
is a right $\infty$-derivator. 
\end{example}

\begin{definition}
An $\infty$-prederivator $\der \colon \dia^{\op} \to \catinf$ is an $\infty$-\emph{derivator} if it is both a left and a right $\infty$-derivator.
\end{definition}

We also specialize these definitions to $n$-prederivators as follows. 

\begin{definition}
An $\infty$-prederivator $\der \colon \dia^{\op} \to \catinf$ is a (left, right) $n$-\emph{derivator} if it is an $n$-prederivator and a (left, right) $\infty$-derivator. 
\end{definition}

\begin{example}
A (left, right) 1-derivator with domain $\mathsf{Cat_1}$ is a (left, right) derivator in the usual sense \cite{Malt, Groth}. Indeed (Der 1)--(Der 3) are completely analogous to the usual axioms. (Der 4) is a convenient variation of the usual axiom for ordinary derivators; the equivalence of the definitions is based on known properties of homotopy exact squares in the context of ordinary derivators (see \cite{Groth, Malt2}). 
\end{example}

\begin{example}
Let $\der \colon \dia^{\op} \to \catn$ be an $n$-prederivator, where $n \in \mathbb{Z}_{\geq 1} \cup \{\infty\}$. For any $k < n$, there is an associated 
$k$-prederivator:
$$\h_k \der \colon \dia^{\op} \to \catn \xrightarrow{\h_k} \mathsf{Cat}_k.$$
If $\der$ is a left (right) $n$-derivator, then $\h_k \der$ is a left (right) $k$-derivator.
\end{example}

As a consequence of the axioms of Definition \ref{deriv_def}, we have the following useful strong version of (Der 4) and (Der 4)*, which identifies a larger 
class of squares for which the base change transformations are equivalences. Similar results are known for $\infty$-categories and for ordinary 
($1$-)derivators (see, for example, \cite{Cisinski} and \cite{Malt2, Groth}).

\begin{proposition} \label{Beck_Chevalley} 
Let $\der$ be an $\infty$-derivator with domain $\dia$. Consider a pullback square in $\dia$:
$$
\xymatrix{
\z \ar[d]_{p}\ar[r]^{j} & \x \ar[d]^-{u}\\
\w \ar[r]_{q} & \y. 
}
$$
\begin{itemize}
\item[(1)] The canonical base change transformation:
$$p_! j^* \longrightarrow q^* u_!$$
is a natural equivalence if $u$ is a cocartesian fibration or if $q$ is a cartesian fibration.
\item[(2)] The canonical base change transformation:
$$q^* u_* \longrightarrow p_* j^*$$
is a natural equivalence if $u$ is a cartesian fibration or if $q$ is a cocartesian fibration.
\end{itemize}
\end{proposition} 
\begin{proof}
Using the duality $\der \mapsto \der(-^{\op})^{\op}$, it suffices to prove only (1). Indeed we obtain (2) by applying (1) to the derivator $\der(-^{\op})^{\op}$ in the case of the opposite pullback square in $\dia$. 

\smallskip

Suppose that $u$ is a cocartesian fibration. Applying (Der 2) and using the naturality properties of base change transformations, it suffices to prove the claim only in the case where $W = \Delta^0$:
$$
\xymatrix{
 \x_y \ar[d]_{p}\ar[r]^{j} & \x \ar[d]^-{u}\\
\Delta^0 \ar[r]_{i_y} & \y. 
}
$$
Using the following factorization of this square and applying (Der 4), it suffices to prove the claim for the left square in the following diagram
$$
\xymatrix{
\x_{y} \ar[r]^j \ar[d]_{p} & \x_{/y} \ar[d]_{p_{u/y}} \ar[r]^{j_{u / y}} & \x \ar[d]^-{u}\\
\Delta^0 \ar[r]_{i_{(y = y)}} & \y_{/y} \ar[r]_{q_{\y / y}} & \y. 
}
$$ 
Note that the horizontal functors in the left square admit left adjoints 
$$\ell_1 \colon \x_{/ y} \to \x_y$$
$$\ell_2 \colon \y_{/y} \to \Delta^0$$
because $u$ is a cocartesian fibration and $\y_{/y}$ has a terminal object $(y = y)$. Thus, after applying $\der$, we obtain pairs of adjoint functors $(j^*, \ell_1^*)$ and $(i_{(y=y)}^*, \ell_2^*)$. Therefore the base change transformation between compositions of left adjoints:  
\begin{equation} \label{base_change_simple}
p_! j^* \longrightarrow i_{(y=y)}^* (p_{u/y})_!
\end{equation}
is conjugate to the natural equivalence
$$\mathrm{id} \colon  p_{u/y}^* (\ell_2)^* \simeq  (\ell_1)^* p^*$$
which then implies that \eqref{base_change_simple} is an equivalence, as required. 

\medskip 

\noindent Suppose now that $q$ is a cartesian fibration. Then it suffices to show that the conjugate base change transformation
\begin{equation} \label{conj_base_change}
u^*q_* \to j_*p^*
\end{equation}
is a natural equivalence. Again, by (Der 2), it suffices to restrict to the case where $X = \Delta^0$:
$$
\xymatrix{
\z = \w_y \ar[r]^(.6)j \ar[d]^p & \Delta^0 \ar[d]^{u = i_y} \\
\w \ar[r]_q & \y.
}
$$
The proof of \eqref{conj_base_change} in this case is obtained similarly by dualizing the arguments in the previous proof. 
\end{proof}

In the context of ordinary derivators, there is a well-known additional axiom (Der 5) which is very useful in practice. Before we state the generalization of this axiom to higher derivators, we define (similarly to the classical case) for every $\infty$-prederivator $\der$ and $\x, \y \in \dia$, the \emph{underlying $(\x-)$diagram functor}
$$\diag_{\x, \y} \colon  \der(\x \times \y) \rightarrow \Fun(\x, \der(\y))$$
to be the adjoint of the following composition:
$$\x \cong \Fun(\Delta^0, \x) \xrightarrow{(- \times \y)} \Fun(\y, \x \times \y) \xrightarrow{\der} \Fun(\der(\x \times \y), \der(\y)).$$
We say that a functor $F\colon \C\to \D$ between $\infty$-categories 
is $n$--\emph{full} if it restricts to $(n-1)$-connected maps between the mapping spaces. For example, a full functor between ordinary categories is $1$-full 
in this sense. 

\begin{definition}
Let $\der$ be an $\infty$-prederivator and let $n \in \mathbb{Z}_{\geq 1} \cup \{\infty\}$. We say that $\der$ is \emph{n-strong} if the following axiom is satisfied:
\begin{enumerate}[labelindent=0pt,  itemindent=1em]
\item[(Der 5$_n$)] For every $\x \in \dia$, the underlying diagram functor 
$$\diag_{\Delta^1, \x} \colon \der(\Delta^1 \times \x) \rightarrow \Fun(\Delta^1, \der(\x))$$
is $n$-full and essentially surjective. 
\end{enumerate}
\end{definition}

\begin{remark} \label{strong-Der5}
Similarly to the case of ordinary (pre)derivators, we may also consider a stronger form of the last axiom (cf. \cite{Heller}) which states that the underlying diagram functor
$$\diag_{\mathcal{I}, \x} \colon \der(\mathcal{I} \times \x) \rightarrow \Fun(\mathcal{I}, \der(\x))$$
is $n$-full and essentially surjective for every $\mathcal{I}$ in $\dia$ which is equivalent in the Joyal model category to a finite 1-dimensional simplicial set. 
Moreover, the following even stronger form of this axiom fits naturally in the higher categorical context:
\begin{enumerate}[labelindent=0pt,  itemindent=1em]
\item[($\widetilde{\text{Der 5}_n}$)] For every $\x \in \dia$ and for every $\mathcal{I} \in \dia$ which is equivalent (in the Joyal model category)
to a $d$-dimensional simplicial set for some $d \leq n + 1$, the underlying diagram functor 
$$\diag_{\mathcal{I}, \x} \colon \der(\mathcal{I} \times \x) \rightarrow \Fun(\mathcal{I}, \der(\x))$$
is $(n+1-d)$-full and essentially surjective. 
\end{enumerate}
Note that this axiom for $d = 0$ is also closely related to (Der 1). 
\end{remark}

Following the definitions of pointed and stable (1-)derivators \cite{Malt, Groth}, we also define pointed and stable $\infty$-(pre)derivators as follows. 

\begin{definition}
An $\infty$-prederivator $\der \colon \dia^{\op} \to \catinf$ is called \emph{pointed} if it lifts to the $\infty$-category $\catinf_{,*}$ of pointed $\infty$-categories and 
functors which preserve zero objects. An $\infty$-derivator $\der \colon \dia^{\op} \to \catinf$ is called \emph{stable} if it is pointed and the associated 1-derivator 
$\h_1 \der$ is stable.
\end{definition} 

\subsection{Examples of $n$-derivators} The examples of $n$-derivators we are mainly interested in are those where the underlying prederivator arises from an 
$\infty$-category as in Example \ref{represented_der}. Let $\dia \subset \catinf$ be a fixed full subcategory satisfying the conditions (Dia 0)--(Dia 3). 

\begin{proposition} \label{represented_preder_prop}
Let $\C$ be an $\infty$-category and let $n \in \mathbb{Z}_{\geq 1} \cup \{\infty \}$. 
\begin{itemize}
\item[(a)] The $n$-prederivator $\der^{(n)}_{\C}$ satisfies (Der 1), (Der 2) and (Der 5$_n$). Moreover, the underlying diagram functor 
$$\diag_{\mathcal{I}, \x} \colon \der^{(n)}_{\C}(\mathcal{I} \times \x) \rightarrow \Fun(\mathcal{I}, \der^{(n)}_{\C}(\x))$$
is $(n + 1 -d)$-full and essentially surjective for every $\mathcal{I}$ in $\dia$ which is equivalent (in the Joyal model category) to a $d$-dimensional simplicial set for some $d \leq n + 1$ (i.e., $\der^{(n)}_{\C}$ satisfies $(\widetilde{\text{Der 5}_n})$, too). 
\item[(b)] Suppose that $\C$ admits $\x$-colimits (resp. $\x$-limits) for any $\x \in \dia$. Then $\der^{(n)}_{\C}$ satisfies (Der 3) and (Der 4) (resp. (Der 3)* and (Der 4)*). As a 
consequence, $\der^{(n)}_{\C}$ is a right (resp. left) $n$-derivator which is $n$-strong. 
\end{itemize}
\end{proposition}
\begin{proof}
(a) (Der 1) is obvious. (Der 2) says that the equivalences in $\h_n(\C^{\x})$, $\x \in \dia$, are given pointwise, which holds by a theorem of Joyal \cite[Chapter 5]{Joyal_1}. 
(Der 5$_n$) follows from Corollary \ref{n-d_homotopy_cat}. The second claim also follows from Corollary \ref{n-d_homotopy_cat} (and Lemma \ref{lem:simplices}) because we may replace $\mathcal{I}$ by a $d$-dimensional 
simplicial set. 

(b) It suffices to show that (Der 3) and (Der 4) hold for the $\infty$-prederivator $\der^{(\infty)}_{\C}$, that is, it suffices to show that $\C$ admits Kan extensions along functors in $\dia$ and that Kan extensions are 
given pointwise by the usual formulas as axiomatized in (Der 4). These claims are established in \cite[6.4]{Cisinski} (see also \cite[4.3.2--4.3.3]{LurHTT}). 
\end{proof}

It is often possible to reduce statements about $\infty$-derivators to the corresponding (known) statements about $1$-derivators. This happens, for example, 
in the case of statements which involve the detection of equivalences. The following proposition shows another instance of this phenomenon and produces 
many examples of $n$-derivators from known examples of $1$-derivators. 

\begin{theorem} \label{deriv_criterion}
Let $\C$ be an $\infty$-category. The following are equivalent:
\begin{itemize}
\item[(1)] $\C$ admits $\x$-colimits and $\x$-limits for every $\x \in \dia$. 
\item[(2)] $\der^{(\infty)}_{\C}$ is an $\infty$-derivator.
\item[(3)] $\der^{n}_{\C}$ is an $n$-derivator for every $n \in \mathbb{Z}_{\geq 1}$.
\item[(4)] $\der^{n}_{\C}$ is an $n$-derivator for some $n \in \mathbb{Z}_{\geq 1}$.
\item[(5)] $\der_{\C} = \der^{(1)}_{\C}$ is an $1$-derivator.
\end{itemize}
\end{theorem}
\begin{proof}
(1) $\Rightarrow$ (2) was shown in Proposition \ref{represented_preder_prop} and the implications (2) $\Rightarrow$ (3) $\Rightarrow$ (4) $\Rightarrow$ (5) are obvious. We prove the implication (5) $\Rightarrow$ (1).  We restrict to showing that $\C$ 
admits $\x$-colimits as the case of $\x$-limits can be treated similarly by duality. Let $F\colon \x \to \C$ 
be an $\x$-diagram and let $G \colon \x^{\triangleright} \to \C$ be the diagram which is the image of $F \in \der_{\C}(\x)$ under 
$$u_! \colon \der_{\C}(\x) \longrightarrow \der_{\C}(\x^{\triangleright})$$ 
where $u \colon \x \to \x^{\triangleright}$ is the canonical inclusion. As a consequence of (Der 4), the functor $u_!$ is fully faithful because $u$ is so. In particular, there is a canonical isomorphism $F \cong u^* u_!(F),$ and therefore we may assume that $G$ is an extension of the diagram $F$. We claim 
that $G$ is a colimit diagram in $\C$ for the functor $F$. For this, it suffices to prove that for each sieve between (nerves of) finite posets  $v \colon \y' \to \y$ in $\dia$, the canonical map 
\begin{equation} \label{comparison_map}
\xymatrix{
\Map(\y, \Map_{\C^{\x^{\triangleright}}}(G, c_y)) \ar[d] \\
\Map(\y', \Map_{\C^{\x^{\triangleright}}}(G, c_y))  \times_{\Map(\y', \Map_{\C^{\x}}(F, c_y))} \Map(\y, \Map_{\C^{\x}}(F, c_y)) 
}
\end{equation}
is a $\pi_0$-isomorphism, where $c_y$ denotes the constant functor at $y \in \C$ in the respective functor $\infty$-category. We will do this by first expressing $\pi_0$ of the domain and 
the target of this map \eqref{comparison_map} in terms of morphism sets in the (1-)category $\der_{\C}(\x^{\triangleright} \times \y)$ and then use the derivator properties of $\der_{\C}$. 

\medskip

\noindent First, note that we have a canonical isomorphism of morphism sets:
\begin{equation} \label{domain_comparison_map}
\begin{split}
\der_{\C}(\x^{\triangleright} \times \y)(\pi_{\x^{\triangleright}, \y}^*(G), c_y) = \pi_0 \big(\Map_{\C^{\x^{\triangleright} \times \y}}(\pi_{\x^{\triangleright}, \y}^*(G), c_y)\big) \\ 
\cong \pi_0 \big(\Map(\y, \Map_{\C^{\x^{\triangleright}}}(G, c_y)) \big)
\end{split}
\end{equation}
where $\pi_{\x^{\triangleright}, \y} \colon \x^{\triangleright} \times \y \to \x$ denotes the projection functor. Similarly, we have canonical 
isomorphisms
\begin{equation} \label{target_comparison_map_2}
\begin{split}
\der_{\C}(\x \times \y)(\pi_{\x, \y}^*(F), c_y) \cong \pi_0 \big(\Map(\y, \Map_{\C^{\x}}(F, c_y)) \big)  \\
\der_{\C}(\x \times \y^{'})(\pi_{\x, \y^{'}}^*(F), c_y) \cong \pi_0 \big(\Map(\y^{'}, \Map_{\C^{\x}}(F, c_y)) \big)
\end{split} 
\end{equation}
where $\pi_{\x, \y}$ and $\pi_{\x, \y'}$ denote again the projection functors.

\medskip 

\noindent Then consider the following pullback diagram in $\dia$:
$$
\xymatrix{
\x \times \y \ar[r]^(.6){\pi_{\x, \y}} \ar[d]_{u \times 1} &  \x \ar[d]^u \\
\x^{\triangleright} \times \y \ar[r]_(.6){\pi_{\x^{\triangleright}, \y}} & \x^{\triangleright}. 
}
$$
Since the horizontal functors are cartesian fibrations, it follows from (Der 4) and Proposition \ref{Beck_Chevalley}(1) that the canonical base change morphism in $\der_{\C}(\x^{\triangleright} \times \y)$:
\begin{equation} \label{eqn_1}
\pi_{\x^{\triangleright}, \y}^*(G) = \pi_{\x^{\triangleright}, \y}^* u_!(F) \stackrel{\cong}{\leftarrow} (u \times 1)_! \pi_{\x, \y}^*(F)
\end{equation}
is an isomorphism. 

\medskip

\noindent Let $i \colon \a = \x^{\triangleright} \times \y' \cup_{\x \times \y'} \x \times \y \subset \x^{\triangleright} \times \y$ denote the full subcategory. 
This is again in $\dia$ by (Dia 1) and (Dia 3) because it can be described as the pullback of a functor $\x^{\triangleright} \times \y \to \Delta^1 \times \Delta^1$ along the `upper corner' inclusion $\Delta^1 \cup_{\Delta^0} \Delta^1 \to \Delta^1 \times \Delta^1$. Note that using \eqref{target_comparison_map_2}, we can identify the set of components of the target of \eqref{comparison_map} canonically with the morphism set:
\begin{equation} \label{target_comparison_map}
\begin{split}
\der_{\C}(\a)(i^* \pi^*_{\x^{\triangleright}, \y}(G), c_y) = \der_{\C}(\a)(i^*\pi_{\x^{\triangleright}, \y}^*u_!(F), c_y) \\
\cong \pi_0(\textsc{target of } \eqref{comparison_map})
\end{split}
\end{equation} 
and the map \eqref{comparison_map} on $\pi_0$ agrees using the identifications \eqref{domain_comparison_map} and \eqref{target_comparison_map} with the map defined by the restriction functor 
$$i^* \colon \der_{\C}(\x^{\triangleright} \times \y) \to \der_{\C}(\a).$$  
Since $i$ is full, it follows from (Der 4) that the unit transformation $1 \to i^* i_!$ of the adjunction $(i_!, i^*)$ is a natural isomorphism.  Therefore 
it suffices to show that the counit morphism 
\begin{equation} \label{eqn_3} 
i_! i^* \big( \pi^*_{\x^{\triangleright}, \y}(G)\big) \to \pi^*_{\x^{\triangleright}, \y} (G)
\end{equation}
is an isomorphism. 

\medskip 

\noindent Consider the following pullback diagram in $\dia$:
$$
\xymatrix{
\x \times \y \ar[d]_{j} \ar[rr]^(.6){\pi_{\x, \y}} && \x \ar[d]^u \\
\a \ar[rr]_{q = \pi_{\x^{\triangleright}, \y} i} && \x^{\triangleright}. 
}
$$
The bottom functor $q$ is a cartesian fibration because it is the composition of cartesian fibrations. Therefore it follows from (Der 4) and Proposition \ref{Beck_Chevalley}(1) that the canonical base change morphism in $\der_{\C}(\a)$:
\begin{equation} \label{eqn_2}
q^*(G) = q^*u_!(F) \stackrel{\cong}{\leftarrow} j_! \pi_{\x, \y}^*(F)
\end{equation}
is an isomorphism. As a consequence of \eqref{eqn_1} and \eqref{eqn_2}, we obtain canonical isomorphisms as follows,
\begin{equation*}
\begin{split} 
\pi_{\x^{\triangleright}, \y}^*(G) = \pi_{\x^{\triangleright}, \y}^* u_!(F) \cong (u \times 1)_! \pi_{\x, \y}^*(F)  \\
\cong i_! j_! \pi_{\x, \y}^*(F) \cong i_! q^* u_!(F) \\ 
= i_! i^* \pi^*_{\x^{\triangleright}, \y}(G).
\end{split}
\end{equation*}
This implies that \eqref{eqn_3} is an isomorphism and therefore, using the adjunction $(i_!, i^*)$ as explained above, it follows that the map \eqref{comparison_map} is a $\pi_0$-isomorphism. 
\end{proof} 

\begin{remark}
We point out a significant simplification of the proof of Theorem \ref{deriv_criterion} in the case where $\dia$ is large enough so that it can detect equivalences of $\infty$-groupoids, i.e., in the case where the following holds: a map of $\infty$-groupoids $X \to X'$ is an equivalence if and only if
$$\pi_0(\Map(\y, X)) \longrightarrow \pi_0(\Map(\y, X'))$$
is an isomorphism for every $\y \in \dia$. This happens, for example, when $\dia$ contains all posets -- and not just the finite ones. Assuming that 
$\dia$ is large enough in this sense, then we may restrict to the case $\y' = \varnothing$ in the proof above, in which case the proof becomes more 
immediate. 
\end{remark}

\section{$K$-theory of higher derivators} \label{sec:K-theory_higher_derivators}

\subsection{Preliminaries} We first recall the $\infty$-categorical version of Waldhausen's $\S_{\bullet}$-construction \cite{Waldhausen, BGT}. 
Let $\C$ be a pointed $\infty$-category which admits finite colimits. For every $n \geq 0$, let $\Ar[n]$ denote the (nerve of the) category 
of morphisms of the poset $[n]$. The $\infty$-category $\S_n \C$ is the full subcategory of $\Fun(\Ar[n], \C)$ spanned by the objects $F \colon \Ar[n] \to \C$ 
such that:
\begin{itemize} 
 \item[(i)] $F(i \to i)$ is a zero object for all $i \in [n]$. 
 \item[(ii)] For every $i \leq j \leq k$, the following diagram in $\C$,
\[
 \xymatrix{
 F(i \to j) \ar[r] \ar[d] & F(i \to k) \ar[d] \\
 F(j \to j) \ar[r] & F(j \to k)
 }
\]
is a pushout. 
\end{itemize}
The construction is clearly functorial in $[n]$, $n \geq 0$, and $\S_{\bullet} \C$ defines a simplicial object of pointed $\infty$-categories, which is functorial in $\C$ with respect to functors which preserve zero objects and pushouts. 

We denote by $\S_{\bullet}^{\simeq} \C$ the associated simplicial object of pointed $\infty$-groupoids, which is obtained by passing pointwise to the maximal $\infty$-subgroupoids of $\S_{\bullet} \C$. For $n \geq 1$, the $\infty$-groupoid $\S_{n}^{\simeq} \C$ is equivalent to $\Map(\Delta^{n-1}, \C)$. Moreover, we have $\S^{\simeq}_{0} \C \simeq \Delta^0$, so we may regard the geometric realization $|\S^{\simeq}_{\bullet} \C|$ as canonically pointed by a zero object in $\C$. The Waldhausen $K$-theory of $\C$ is defined to be the space:
$$K(\C) : = \Omega | S_{\bullet}^{\simeq} \C |.$$
If $\C$ arises from a nice Waldhausen category, this definition of $K$-theory agrees up to homotopy equivalence with the Waldhausen $K$-theory of the corresponding Waldhausen category (see \cite{BGT}). The definition of $K$-theory is functorial with respect to functors $F \colon \C \to \C'$ which preserve the zero object and finite colimits. 

\medskip

Following \cite[Lemma 1.4.1]{Waldhausen} and \cite[Proposition 4.2.1]{MR2}, we consider also the following simpler model for Waldhausen $K$-theory. Restricting pointwise to the objects of $\S_{\bullet} \C$, we obtain a simplicial set 
$$\s_{\bullet} \C \colon \Delta^{\op} \to \mathsf{Set}, \ \ [n] \mapsto \s_n \C : = (\S_n \C)_0.$$
There is a canonical comparison map, given by the inclusion of objects,
$$\iota \colon \Omega |\s_{\bullet} \C| \longrightarrow \Omega | \S_{\bullet}^{\simeq} \C| = K(\C).$$

\begin{proposition} \label{comparison_s_constr}
The comparison map $\iota$ is a weak equivalence. 
\end{proposition}
\begin{proof}
The proof is essentially the same as the proof of \cite[Lemma 1.4.1 and Corollary]{Waldhausen} (see also \cite[Proposition 4.2.1]{MR2}).
\end{proof}

\subsection{Derivator $K$-theory} We extend the definition of derivator $K$-theory of Maltsiniotis \cite{Malt} and Garkusha \cite{Garkusha2, Garkusha1} to general pointed  right $\infty$-derivators. As in the case of ordinary derivators, this definition is based on the following intrinsic notion of cocartesian 
square.

\medskip

Let $i \colon \ulcorner = \Delta^1 \cup_{\Delta^0} \Delta^1 \to \square = \Delta^1 \times \Delta^1$ denote the `upper corner' inclusion. For any right $\infty$-derivator $\der$ (with domain $\dia$), we have an adjunction:
$$i_! \colon \der(\ulcorner) \rightleftarrows \der(\square) \colon i^*.$$

\begin{definition}
Let $\der \colon \dia^{\op} \to \catinf$ be a right $\infty$-derivator with domain $\dia$. An object $F \in \der(\square)$ is called \emph{cocartesian} if the canonical morphism 
$$i_! i^*(F) \to F$$
is an equivalence in $\der(\square)$. 
\end{definition}

Let $\der$ be a pointed right $\infty$-derivator (with domain $\dia$). Before we define the $K$-theory of $\der$, we first need to introduce some more notation. For every $0 \leq i \leq j \leq k \leq n$, we denote by $i_{i,j,k} \colon \square 
\to \Ar[n]$ the inclusion of the following square in $\Ar[n]$: 
\[
\xymatrix{
(i \to j) \ar[r] \ar[d] & (i \to k) \ar[d] \\
(j \to j) \ar[r] & (j \to k).
}
\]
We define $\S_n \der$ to be the full subcategory of $\der(\Ar[n])$ which is spanned by the objects $F \in \der(\Ar[n])$ such that:
\begin{itemize} 
\item[(i)] $F_{(i \to i)}$ is a zero object for all $i \in [n]$. 
 \item[(ii)] For every $i \leq j \leq k$, the object $i_{i,j,k}^*(F) \in \der(\square)$, which may be depicted as follows:
\[
 \xymatrix{
 F_{(i \to j)} \ar[r] \ar[d] & F_{(i \to k)} \ar[d] \\
 F_{(j \to j)} \ar[r] & F_{(j \to k)},
 }
\]
is cocartesian in $\der(\square)$. 
\end{itemize}
The assignment $[n] \mapsto \S_n \der$ defines a simplicial object of pointed $\infty$-categories. Moreover, it is natural with respect to strict morphisms between pointed right $\infty$-derivators which preserve zero objects and cocartesian squares. 

Let $\S_{\bullet}^{\simeq} \der$ denote the simplicial object of pointed $\infty$-groupoids,  
which is obtained by passing pointwise to the maximal $\infty$-subgroupoids of $\S_{\bullet} \der$. We have $\S_0^{\simeq} \der \simeq \Delta^0$  and we regard the geometric realization $|\S_{\bullet}^{\simeq} \der|$ as based at a zero object of $\der(\Delta^0)$. 

\begin{definition} \label{derivator_K_theory_def}
Let $\der \colon \dia^{\op} \to \catinf$ be a pointed right $\infty$-derivator with domain $\dia$. The \emph{derivator $K$-theory} of $\der$ is defined to be the space: 
$$K(\der) : = \Omega | \S_{\bullet}^{\simeq} \der|.$$
\end{definition} 

We note that the definition of derivator $K$-theory is functorial with respect to strict morphisms $F \colon \der \to \der'$ which preserve the zero object and cocartesian squares. Moreover, derivator $K$-theory is invariant under those strict morphisms which are pointwise equivalences of $\infty$-categories. 

This definition of derivator $K$-theory clearly agrees with the usual derivator $K$-theory for ordinary derivators \cite{Garkusha2, MR2}. Moreover, this definition of derivator $K$-theory is also an \emph{extension} of the $K$-theory of pointed $\infty$-categories with finite colimits to pointed right $\infty$-derivators; by definition, for any pointed $\infty$-category $\C$ with finite colimits, the $K$-theory of $\C$ agrees with the derivator $K$-theory of $\der^{(\infty)}_{\C}$.  

\begin{remark} Waldhausen's Additivity Theorem \cite{Waldhausen} establishes one of the fundamental properties of Waldhausen $K$-theory. The analogue of this 
theorem has been established for the derivator $K$-theory of stable $1$-derivators by Cisinski--Neeman \cite{C-N}, confirming one of Maltsiniotis' 
conjectures in \cite{Malt}. This theorem has recently been extended to general pointed right 1-derivators by Coley \cite{Coley}. It would be interesting to know if the additivity theorem holds more generally for the derivator $K$-theory of pointed right $\infty$-derivators. 
\end{remark}

\subsection{Comparison with Waldhausen $K$-theory} \label{mu-n-map}
Let $\C$ be a pointed $\infty$-category with finite colimits. Applying pointwise the homotopy $n$-category functor  to the simplicial object $[k] \mapsto \S^{\simeq}_k \C$, we obtain a new simplicial object of (pointed) $\infty$-groupoids, 
$$\h_{n} \S_{\bullet}^{\simeq} \C: \Delta^{\op} \to \mathsf{Grpd}_{\infty}, \ [k] \mapsto \h_{n} \big(\S^{\simeq}_k \C\big),$$
and there is a canonical comparison map:
$$\S_{\bullet}^{\simeq} \C \longrightarrow \h_n \big(\S_{\bullet}^{\simeq} \C\big).$$
(We note that the $n$-groupoid $\h_n(\D^{\simeq})$ is the maximal $\infty$-subgroupoid of $\h_n(\D)$ for every $\infty$-category $\D$; this can be seen directly from the definitions of $\h_n(-)$ and $(-)^{\simeq}$. As a result, the order of the operations $\h_n(-)$ and $(-)^{\simeq}$ does not play an essential role in the definition of $\h_n \S_{\bullet}^{\simeq} \C$.) 

Let $\der^{(n)}_{\C}$ be the pointed right $n$-derivator associated to $\C$ with domain $\dirf$ (see Proposition \ref{represented_preder_prop}). As a consequence of the natural identification between $\h_n \big(\S_{\bullet}^{\simeq} \C \big)$ and $\S_{\bullet}^{\simeq} \der^{(n)}_{\C},$
we obtain a canonical comparison map from the Waldhausen $K$-theory of $\C$ to derivator $K$-theory:
$$\mu_n: K(\C) \to K(\der^{(n)}_{\C}).$$
In addition, the natural morphisms of simplicial objects $\h_n\big(\S_{\bullet}^{\simeq} \C\big) \to \h_{n-1}\big(\S_{\bullet}^{\simeq} \C\big)$, for $n > 1$, define a tower of derivator $K$-theories for $\C$ which is compatible with the comparison maps $\mu_n$:
$$
\xymatrix{
& K(\C) \ar[d]_{\mu_n} \ar[dr]_{\mu_{n-1}} \ar[drrrr]^{\mu_1} &&&& \\
\cdots \ar[r] & K(\der^{(n)}_{\C}) \ar[r] &  K(\der^{(n-1)}_{\C}) \ar[r] & \cdots \ar[rr] && K(\der^{(1)}_{\C}).
}
$$
In the case of ordinary derivators, Maltsiniotis conjectured in \cite{Malt} that the comparison map $\mu_1$ is a weak equivalence for exact categories. The comparison map $\mu_1$ was subsequently studied in \cite{Garkusha1, Muro, MR2, MR1}. It is known (see \cite{MR1}) that $\mu_1$ is not a weak equivalence for general $\C$, and moreover, that $\mu_1$ will fail to be a weak equivalence even for exact categories \emph{if} derivator $K$-theory satisfies \emph{localization} -- a property which was also conjectured by Maltsiniotis \cite{Malt}. 

We prove a general result about the connectivity of the comparison map $\mu_n$ for general $\C$ and $n \geq 1$. This connectivity estimate is also a small improvement of the known estimate for $n=1$ that was shown by Muro \cite{Muro}. 

\begin{theorem} \label{comparison_K-theory}
Let $\C$ be a pointed $\infty$-category which admits finite colimits. Then the comparison map $\mu_n: K(\C) \to K(\der^{(n)}_{\C})$
 is $(n+1)$-connected. 
\end{theorem}

We will need the following useful elementary fact about simplicial spaces. 

\begin{lemma} \label{technical_lemma}
Let $f_{\bullet} \colon X_{\bullet} \to Y_{\bullet}$ be a map of simplicial spaces. Suppose that $f_k$ is $(m-k)$-connected for every $k \geq 0$. Then the 
map $||f_{\bullet}|| \colon ||X_{\bullet}|| \to ||Y_{\bullet}||$ is $m$-connected. (Here $||-||$ denotes the fat geometric realization of a simplicial space.) 
\end{lemma}
\begin{proof}
See \cite[Lemma 2.4]{E-RW}.
\end{proof}

\begin{proof}(of Theorem \ref{comparison_K-theory})
By Lemma \ref{technical_lemma}, it suffices to show that the map of $\infty$-groupoids
$$\S_k^{\simeq} \C \to \h_{n} \big(\S_k^{\simeq} \C\big)$$
is $(n+2-k)$-connected for all $k \geq 0$. This holds since the map is an equivalence for $k = 0$ and $(n+1)$-connected for $k > 0$ (Example \ref{truncation_groupoid}).
\end{proof}

\begin{remark} \label{best_possible?}
The main result of \cite[Theorem 1.2]{MR1} shows that the comparison map $\mu_1$ is not a $\pi_3$-isomorphism in general. (In addition, a closer inspection of the proofs in \cite{MR1} also shows that the map $\mu_1$ will not be 3-connected if derivator $K$-theory satisfies localization.)  
It seems likely that the connectivity estimate of Theorem \ref{comparison_K-theory} is best possible in general. 
\end{remark}

\begin{remark} \label{t-structure}
By \cite[Theorem 7.1]{Garkusha1}, the comparison map $\mu_1$ is $\pi_*$-split injective in the case where $\C$ is the bounded derived category of an abelian category. In fact, it is shown \cite{Garkusha1} that there is 
a retraction map to $\mu_1$. As a consequence, the comparison map $\mu_n$ also admits a retraction in this case for all $n \geq 1$. Related to this, an interesting problem suggested by B. Antieau is whether $\mu_n$ is a weak equivalence when $\C$ is a stable $\infty$-category which admits a bounded t-structure. 
\end{remark}

\subsection{Waldhausen $K$-theory of derivators} 

Waldhausen $K$-theory for pointed right 1-derivators was defined in \cite{MR2}. It was shown in \cite{MR2} that this $K$-theory of derivators agrees with the usual Waldhausen $K$-theory for all well-behaved Waldhausen categories \cite[Theorem 4.3.1]{MR2}. We consider an analogous definition of $K$-theory for general pointed right $\infty$-derivators.

\medskip 

Let $\der$ be a pointed right $\infty$-derivator with domain $\dia$.  Let $\S_{\bullet \bullet} \der$ be the bisimplicial set whose set of $(n,m)$-simplices $\S_{n,m} \der$ is the set of objects 
$$F \in \ob \big(\der(\Delta^m \times \Ar[n])\big)$$
 such that:
\begin{itemize}
\item[(1)] for every $j\colon [0] \to [m]$ the object $(j \times\mathrm{id})^*(F) \in \Ob \big(\der(\Ar[n])\big)$ is in $\S_n \der$.
\item[(2)] the underlying diagram functor associated to $F$,
$$\diag_{\Delta^m, \Ar[n]}(F) \colon \Delta^m \to \der(\Ar[n]),$$
takes values in equivalences. 
\end{itemize}
The bisimplicial operators of $\S_{\bullet \bullet} \der$ are again defined using the structure of the underlying $\infty$-prederivator. Moreover, it is easy to see that the construction is functorial in $\der$ with respect to strict morphisms which preserve the zero objects and cocartesian squares. 
We regard the geometric realization $| \S_{\bullet \bullet} \der|$ as based at a zero object of $\der(\Delta^0)$. 
 
\begin{definition}
Let $\der \colon \dia^{\op} \to \catinf$ be a pointed right $\infty$-derivator with domain $\dia$. The \emph{Waldhausen $K$-theory} of $\der$ is defined to be the space:
$$K^W(\der) := \Omega | \S_{\bullet \bullet} \der|.$$
\end{definition}

Following \cite{MR2}, we consider also the analogue of the $\s_{\bullet}$-construction in this context. Restricting pointwise to the objects of $\S_{\bullet} \der$, we obtain a simplicial set 
$$\s_{\bullet} \der \colon \Delta^{\op} \to \mathsf{Set}, \ \ [n] \mapsto \s_n \der := \S_{n, 0} \der = (\S_n \der)_0.$$
There is a canonical comparison map, given by the inclusion of $0$-simplices,
$$\iota \colon \Omega |\s_{\bullet} \der| \longrightarrow \Omega | \S_{\bullet \bullet} \der| = K^W(\der).$$
This is the analogue of the comparison map in Proposition \ref{comparison_s_constr} for pointed right $\infty$-derivators.

\begin{proposition} \label{Waldhausen_K_theory_der} (a) The comparison map $\iota$ is a weak equivalence. (b)
Let $\C$ be a pointed $\infty$-category with finite colimits and let $n \in \mathbb{Z}_{\geq 1} \cup \{\infty\}$. There is a commutative diagram of weak equivalences:
$$
\xymatrix{
\Omega |\s_{\bullet} \C| \ar[d]_{\simeq} \ar@{=}[rr] && \Omega |\ob \S_{\bullet} \der^{(n)}_{\C} | = \Omega | \s_{\bullet} \der^{(n)}_{\C} | \ar[d]^{\simeq}_{\iota} \\
K(\C) = \Omega |\S_{\bullet}^{\simeq} \C| \ar[rr]_{\simeq} && \Omega |\S_{\bullet \bullet} \der^{(n)}_{\C} | = K^W(\der^{(n)}_{\C}). 
}
$$
\end{proposition}
\begin{proof}
(a) The proof is essentially the same as the proof of \cite[Proposition 4.2.1]{MR2} (see also \cite[Lemma 1.4.1 and Corollary]{Waldhausen}). (b) The bottom map is a weak equivalence (independently of $n$!) because we have $(\S_k^{\simeq} \C)_m = \S_{k,m} \der^{(n)}_{\C}$. The left vertical map is a weak equivalence 
by Proposition \ref{comparison_s_constr}. Then the result follows. 
\end{proof}

\subsection{Universal property of derivator $K$-theory}
The comparison maps $\{\mu_n\}$ from Waldhausen $K$-theory to derivator $K$-theory can be defined more generally for pointed right $\infty$-derivators. Given a pointed right $\infty$-derivator $\der$ (with domain $\dia$), 
the underlying diagram functors define a bisimplicial map as follows,
$$\diag_{\Delta^m, \Ar[n]} \colon \S_{n,m} \der \to (\S_{n}^{\simeq} \der)_m,$$
which after passing to the geometric realization and taking loop spaces defines a comparison map:
$$\mu \colon \Omega |\s_{\bullet} \der| \stackrel{(\ref{Waldhausen_K_theory_der})}{\simeq} K^W(\der) \to K(\der).$$
A universal property of this comparison map in the case of $1$-derivators was shown in \cite[Theorem 5.2.2]{MR2}. More specifically, it was shown that $\mu$ is homotopically initial among all natural transformations from Waldhausen $K$-theory to a 
functor which is invariant under (pointwise) equivalences of pointed right derivators. The proof of this universal property extends similarly to our present $\infty$-categorical context. 

\medskip

Let $\Der$ denote the (ordinary) category of pointed right $\infty$-derivators and strict morphisms which preserve the zero object and cocartesian squares. It will be convenient to work with the simpler model for the Waldhausen $K$-theory of derivators given by the $\s_{\bullet}$-construction. This will be denoted by
$$K^{W, \ob} \colon \Der \to \mathrm{Top}, \ \ \der \mapsto \Omega | \s_{\bullet} \der|,$$
where $\mathrm{Top}$ denotes the ordinary category of topological spaces. Then we may regard the comparison map $\mu$ as a natural transformation $K^{W, \ob} \to K$ between functors defined on $\Der$.

\begin{definition} \label{facto_def}
The \emph{category $\facto$ of invariant approximations} to Waldhausen $K$-theory is the full subcategory of the comma category $K^{W, \ob} \downarrow \mathrm{Top}^{\Der}$ 
spanned by the objects $(\eta \colon K^{W, \ob} \to F)$ such that $F \colon  \Der \to \mathrm{Top}$ sends pointwise equivalences in $\Der$ to weak equivalences. A morphism in 
$\facto$ 
$$
\xymatrix{
& K^{W, \ob} \ar[ld]_\eta \ar[rd]^{\eta'} &\\
F\ar[rr]_u&& F'
}
$$
is a \emph{weak equivalence} if the components of $u$ are weak equivalences.
\end{definition} 

\begin{remark}
As defined in Definition \ref{facto_def}, the category $\facto$ (denoted by $\mathtt{App}$ in \cite{MR2}) may not be locally small. This set--theoretical issue can be addressed by restricting to suitable small subcategories of $\Der$, as was done in \cite{MR2}. 
\end{remark}

We recall from \cite{DHKS} that an object $x \in \mathcal{C}$ in an (ordinary) category with weak equivalences $(\mathcal{C}, \mathcal{W})$ (satisfying in addition the ``2-out-of-6'' property) is \emph{homotopically initial} if there are functors $F_0, F_1 \colon \mathcal{C} \to \mathcal{C}$ which preserve the weak equivalences and a natural transformation $\phi \colon F_0 \Rightarrow F_1$ such that: 
\begin{itemize} 
 \item[(i)] $F_0$ is naturally weakly equivalent to the constant functor at $x \in \mathcal{C}$. 
 \item[(ii)] $F_1$ is naturally weakly equivalent to the identity functor on $\mathcal{C}$. 
 \item[(iii)] $\phi_x \colon F_0(x) \to F_1(x)$ is a weak equivalence.
\end{itemize}
A homotopically initial object in $(\mathcal{C}, \mathcal{W})$ defines an initial object in the associated $\infty$-category. 

\begin{theorem}\label{universal_prop_derivator_K-theory}
The object $(\mu \colon K^{W, \ob} \to K) \in \facto$ is homotopically initial. 
\end{theorem}
\begin{proof}(Sketch)
The proof is similar to \cite[Theorem 5.2.2]{MR2} so we only give a sketch of the proof. Given $\der \in \Der$ and $m \geq 0$, let $\der^{\simeq}_{m}$ denote the $\infty$-prederivator whose value at $\x \in \dia$ 
is the full subcategory $\Fun_{\simeq}(\Delta^m, \der(\x)) \subset \Fun(\Delta^m, \der(\x))$ spanned by the functors $\Delta^m \to \der(\x)$ which take values in equivalences. This $\infty$-prederivator is pointwise equivalent
to $\der$ and therefore also a pointed right $\infty$-derivator. Varying $m \geq 0$, we obtain a simplicial object 
$(\der^{\simeq}_{m})_{m \geq 0}$ in $\Der$ with $\der^{\simeq}_0 = \der$.

\medskip 

\noindent For the proof of the theorem, it suffices to note that every object in $\facto$,
$$(\eta_{\der} \colon K^{W, \ob}(\der) \to F(\der))_{\der \in \Der}$$ 
is naturally weakly equivalent (as object in $\facto$) to the composite functor
$$\big(K^{W, \ob}(\der) \to ||K^{W, \ob}(\der^{\simeq}_{\bullet})|| \xrightarrow{||\eta||} ||F(\der^{\simeq}_{\bullet})||\big)_{\der \in \Der}.$$
This uses the fact that $F$ respects pointwise equivalences in $\Der$. Moreover, the first map above defines a natural transformation which is canonically identified with $\mu$. As a result, we have constructed a zigzag of natural transformations 
from the constant endofunctor at $\mu$ to $\mathrm{id}_{\facto}$ satisfying (i)--(iii).
\end{proof}

\section{$K$-theory of homotopy $n$-categories} \label{sec:K-theory_higher_homotopy_cats}

\subsection{Revisiting the properties of homotopy $n$-categories} \label{stable-n-categories} Let $\C$ be a pointed $\infty$-category with finite colimits. Then the associated homotopy $n$-category $\h_n \C$ satisfies the following:
\begin{itemize}
\item[(a)] $\h_n \C$ is a pointed $n$-category.
\item[(b)] The suspension functor $\Sigma_{\C} \colon \C \to  \C$ induces a functor $\Sigma \colon \h_n \C \to \h_n \C$. This is an equivalence if and only if $\C$ is stable (see \cite[Corollary 1.4.2.27]{LurHA}).
\item[(c)] $\h_n \C$ admits finite coproducts and weak pushouts of order $n-1$. Moreover, these are preserved by the functor $\gamma_n \colon \C \to \h_n \C$ (Proposition \ref{preservation of weak colimits}). 
\item[(d)] For every $x \in \h_n \C$, there is a natural weak pushout of order $n-1$,
$$
\xymatrix{
x \ar[d] \ar[r] & 0 \ar[d] \\
0 \ar[r] & \Sigma x.
}
$$
\end{itemize}
Assuming that $\C$ is a stable $\infty$-category, then the adjoint equivalence $(\Sigma_{\C}, \Omega_{\C})$ induces an adjoint equivalence $\Sigma \colon \h_n \C \rightleftarrows \h_n \C \colon \Omega$. Moreover, by duality, $\h_n \C$ also satisfies in this case the following dual versions of (c)--(d):
\begin{itemize}
\item[(c)$'$] $\h_n \C$ admits finite products and weak pullbacks of order $n-1$. Moreover, these are preserved by the functor $\gamma_n \colon \C \to \h_n \C$. 
\item[(d)$'$] For every $x \in \h_n \C$, there is a natural weak pullback of order $n-1$
$$
\xymatrix{
\Omega x \ar[d] \ar[r] & 0 \ar[d] \\
0 \ar[r] & x.
}
$$
\end{itemize}
In addition, if $\C$ is a stable $\infty$-category, $\h_n \C$ satisfies the following property:
\begin{itemize}
\item[(e)] A square in $\h_n \C$ is a weak pushout of order $n-1$ if and only if it is a weak pullback of order $n-1$. 
\end{itemize}
Note that if $n > 1$, weak pushouts (resp. weak pullbacks) of order $n-1$ are unique up to (non-canonical) equivalence. This observation can be used to deduce (e) for $n > 1$. The validity of (e) for $n=1$ can be verified by a direct argument. 
 
\subsection{Towards stable $n$-categories} \label{stable-n-categories2} An attempt towards an axiomatization of the properties (a)--(d) would naturally lead to considering triples
$$(\D, \Sigma \colon \D \to \D, \sigma \colon \D \to \D^{\square})$$ 
where:
\begin{itemize} 
\item[(1)] $\D$ is a pointed $n$-category,
\item[(2)] $\Sigma \colon \D \to \D$ is an endofunctor,
\item[(3)] $\D$ admits finite coproducts and weak pushouts of order $n-1$, 
\item[(4)] The functor $\sigma$ sends an object $x \in \D$ to a weak pushout of order $n-1$ which has the following form:
$$
\xymatrix{
x \ar[d] \ar[r] & 0 \ar[d] \\
0 \ar[r] & \Sigma x.
}
$$
\end{itemize}
Then, specializing to the stable context and using (c)$'$--(d)$'$ and (e), it would be natural to require in addition:
\begin{itemize} 
\item[(2)$'$] $\Sigma \colon \D \to \D$ is an equivalence,
\item[(3)$'$] $\D$ admits finite products and weak pullbacks of order $n-1$, 
\item[(5)  ] A square in $\D$ is a weak pushout of order $n-1$ if and only if it is a weak pullback of order $n-1$. 
\end{itemize}

\medskip

\noindent Antieau \cite[Conjecture 8.28]{Antieau} conjectured that there is a good theory of \emph{stable $n$-categories}, $1 \leq n \leq \infty$, satisfying the following properties:
\begin{itemize}
 \item[(i)] Stable $n$-categories, exact functors and natural transformations form an $(n,2)$-category with a forgetful functor to $\catn$.
 \item[(ii)] For each $n \geq k$, the homotopy $k$-category functor defines a functor from stable $n$-categories to stable $k$-categories.
 \item[(iii)] For $n = \infty$, the theory recovers the theory of stable $\infty$-categories, exact functors, and natural transformations. 
 \item[(iv)] For $n = 1$, the theory recovers the theory of triangulated categories, exact functors, and natural transformations.
\end{itemize}
It seems reasonable to take properties (1)--(5) (incl. (2)$'$--(3)$'$) as a minimal basis for such a notion of stable $n$-category. Firstly, these properties are preserved after passing to lower homotopy categories
(Proposition \ref{preservation of weak colimits}). Also, for any $n > 2$ and any pointed $n$-category $\D$ which satisfies these properties, the 
associated homotopy category $\h_1 \D$ can be equipped with a canonical triangulated structure. The proof is essentially exactly the same as that 
for stable $\infty$-categories in \cite[1.1.2]{LurHA}, using the properties of weak pushouts of order $n-1 > 1$ instead of actual pushouts; we believe that this claim holds also for $n=2$. Also, it would be interesting to explore possible connections between properties of type (1)--(5) and $n$-angulated structures (on the ordinary homotopy category) in the sense of \cite{GKO}. Finally, for $n=\infty$, these properties characterize stable $\infty$-categories. 

\medskip 

On the other hand, concerning the case  $n = 1$, the notion of a triangulated structure includes more structure than what is required in (1)--(5). This could be regarded as a singularity that arises at the lowest level of coherence: since weak pushouts (of order 0) are not unique up to equivalence, they do not yield canonical connecting ``boundary'' maps, not even up to homotopy, and therefore they do not suffice for defining distinguished triangles. Alternatively, it may in fact be desirable to consider additional structure to (1)--(5), in the form of fixed choices of higher weak colimits, and stipulate their properties in analogy with triangulated structures. We will not attempt to give an axiomatic definition of a stable $n$-category in this paper, but we will propose to consider a triple satisfying (1)--(5) as a basic invariant of any good notion of a stable $n$-category.

\subsection{$K$-theory of pointed $n$-categories with distinguished squares} \label{K-theory_distinguished_squares} We will define $K$-theory for certain $n$-categories equipped with distinguished squares. These distinguished squares are meant to play the role of the pushout squares in the definition of Waldhausen $K$-theory. In the case of ordinary categories, this notion of a category with distinguished squares and its $K$-theory can be viewed as a simpler and more basic version of Neeman's $K$-theory of a category with squares as defined in \cite[Sections 5--7]{Neeman}. 

The main example we are interested in is the homotopy $n$-category $\h_n \C$ of a pointed $\infty$-category $\C$ which admits finite colimits, equipped with the squares which come from pushout squares in $\C$ as the distinguished squares. The main motivation for introducing $K$-theory for $\h_n \C$ is to identify a part of Waldhausen $K$-theory $K(\C)$ which may canonically be recovered from $\h_n \C$, regarded as an $n$-category with distinguished squares. In particular, our main result (Theorem \ref{comparison_K-theory_ncat}) generalizes the well-known fact that $K_0(\C)$ can be recovered from $\h_1 \C$, regarded as a category equipped with those squares which arise from pushouts in $\C$.

\begin{definition}
A \emph{pointed} $n$-\emph{category with distinguished squares}, $n \geq 1$, is a pair $(\C, \T)$ where $\C$ is a pointed $n$-category and $\T$ is a collection of weak pushout squares in $\C$ of order $n-1$ which contains the constant squares at a zero object. We call the diagrams in $\T$ \emph{distinguished squares}. 

An \emph{exact functor} $F \colon (\C, \T) \to (\C', \T')$ between pointed $n$-categories with distinguished squares is a functor $F \colon \C \to \C'$ which preserves zero objects and distinguished squares. 
\end{definition}

\begin{definition} \label{canonical-structure}
Let $\C$ be a pointed $\infty$-category with finite colimits and let $n \geq 1$. We define the \emph{canonical structure of distinguished squares in $\h_n \C$} to be the collection of  squares in $\h_n \C$ which are equivalent in $\h_n \C$ to the image of a pushout square in $\C$ under the functor $\gamma_n \colon \C \to \h_n \C$. By Proposition \ref{preservation of weak colimits}, these are weak pushout squares of order $n-1$. For $n>1$, they are precisely the weak pushouts of order $n-1$. We will denote this pointed $n$-category with distinguished squares by $(\h_n \C, \can)$. 
\end{definition}

\begin{remark} \label{canonical-structure-2}
The canonical structure on $\h_n \C$ for $n > 1$ is described in terms of an intrinsic property of $\h_n \C$ and therefore depends only on $\h_n \C$. This 
canonical structure can be defined more generally for any finitely weakly cocomplete $n$-category for $n > 1$ (see Definition \ref{weakly-cocomplete-n-cat}). 
For $n = 1$, the canonical structure is an additional structure on $\h_1 \C$ that is canonically induced from $\C$.  
\end{remark}

Let $(\C,  \T)$ be a pointed $n$-category with distinguished squares. Let $\S_q (\C, \T)$ denote the full subcategory of $\Fun(\Ar[q], \C)$ which is spanned by the functors $F \colon \Ar[q] \to \C$ such that:
\begin{itemize} 
 \item[(i)] $F(i \to i)$ is a zero object for all $i \in [q]$. 
 \item[(ii)] For every $0 \leq i \leq j \leq k \leq m \leq q$, the following diagram in $\C$,
\[
 \xymatrix{
 F(i \to k) \ar[r] \ar[d] & F(i \to m) \ar[d] \\
 F(j \to k) \ar[r] & F(j \to m)
 }
\]
is a distinguished square in $(\C, \T)$. 
\end{itemize}
Note that the construction is functorial in $q \geq 0$ and therefore $\S_{\bullet} (\C, \T)$ defines a simplicial object of pointed $n$-categories. (Of course the construction applies similarly to more general choices $\T$ of distinguished squares, but we will be interested in collections $\T$ which consist of higher weak pushouts.) Moreover, the construction is functorial with respect to exact functors between pointed $n$-categories with distinguished squares. 

We denote by $\S_{\bullet}^{\simeq} (\C, \T)$ the associated simplicial object of $\infty$-groupoids which is obtained by passing to the maximal $\infty$-subgroupoids pointwise.  We have 
$\S_0^{\simeq} (\C, \T) \simeq \Delta^0$  and therefore we may regard the geometric realization $|\S_{\bullet}^{\simeq} (\C, \T)|$ as based
at a zero object of $\C$. 

\begin{definition}
Let $(\C, \T)$ be a pointed $n$-category with distinguished squares. The $K$-\emph{theory of $(\C, \T)$} is defined to be the space: 
$$K(\C, \T) : = \Omega | S_{\bullet}^{\simeq} (\C, \T) |.$$
\end{definition}

\subsection{Comparison with Waldhausen $K$-theory.}
Let $\C$ be a pointed $\infty$-category which admits finite colimits. By Proposition \ref{preservation of weak colimits}, it follows that passing pointwise from $\S_{\bullet} \C$ to the respective homotopy $n$-categories defines map of simplicial objects, 
$$\S_{\bullet} \C \to \S_{\bullet} (\h_n \C, \can),$$
and therefore also a comparison map between $K$-theory spaces
$$\rho_n \colon K(\C) \to K(\h_n \C, \can).$$
Note that this comparison map factors canonically through the comparison map $\mu_n \colon K(\C) \to K(\der^{(n)}_{\C})$ (Subsection \ref{mu-n-map}). 

By Proposition \ref{preservation of weak colimits}, the canonical functors $\gamma_{n-1} \colon \h_n \C \to \h_{n-1} \C$, $n > 1$, define exact functors $(\h_n \C, \can) \to (\h_{n-1} \C, \can)$. Therefore, we obtain a tower of $K$-theories associated to $\C$, which is compatible with the comparison 
maps $\rho_n$,
$$
\xymatrix{
& K(\C) \ar[d]_{\rho_n} \ar[dr]_{\rho_{n-1}} \ar[drrrr]^{\rho_1} &&&& \\
\cdots \ar[r] & K(\h_n \C, \can) \ar[r] &  K(\h_{n-1} \C, \can) \ar[r] & \cdots \ar[rr] && K(\h_1\C, \can).
}
$$
The next result gives a general connectivity estimate for the comparison map $\rho_n$. 

\begin{theorem} \label{comparison_K-theory_ncat}
Let $\C$ be a pointed $\infty$-category which admits finite colimits. For $n \geq 1$, the comparison map $\rho_n \colon K(\C) \to K(\h_n \C, \can)$ 
is $n$-connected. 
\end{theorem}
\begin{proof}
We write $\h_n \C$ for $(\h_n \C, \can)$ and $\S_{\bullet} \h_n\C$ for $\S_{\bullet}(\h_n \C, \can)$, when the canonical structure is understood from the context, in order to simplify the notation in the proof. By Lemma \ref{technical_lemma}, it suffices to show that the map of $\infty$-groupoids
\begin{equation} \label{comparison_map_proof}
\S_q^{\simeq} \C \to \S_q^{\simeq} \h_n \C
\end{equation}
is $(n+1-q)$-connected for every $q \geq 0$. The claim is obvious for $q = 0$. For $q = 1$, the map is $(n+1)$-connected (Example \ref{truncation_groupoid}). For $q = 2$ and $n=1$, the map is $0$-connected by definition. This completes the proof for $n = 1$. 

We may now restrict to the case $n > 1$. We claim that for every $q > 1$, the map \eqref{comparison_map_proof} is $(n-1)$-connected. We consider the diagram
\begin{equation} \label{diagram1}
\xymatrix{
\S_q^{\simeq} \C \ar[rr] \ar[d]_{\simeq} && \S_q^{\simeq} \h_n\C \ar[d] \\
(\C^{\Delta^{q-1}})^{\simeq} \ar[r] & \h_n(\C^{\Delta^{q-1}})^{\simeq} \ar[r] & \big((\h_n \C)^{\Delta^{q-1}}\big)^{\simeq}
}
\end{equation}
where the vertical maps are given by restriction along the inclusion map of posets $[q-1] \subseteq \Ar[q], \ j \mapsto (0 \to j +1).$ The left vertical map in \eqref{diagram1} is an equivalence. The lower left map is $(n+1)$-connected (Example \ref{truncation_groupoid}). The lower right map in \eqref{diagram1} is $n$-connected by Corollary \ref{n-d_homotopy_cat}, since we may replace $\Delta^{q-1}$ by its spine which is $1$-dimensional.

Thus, it suffices to show that the right vertical map in \eqref{diagram1} is $n$-connected for any $q > 0$. This claim is obvious for $q = 1$. For $q > 1$, we proceed by induction on $q$ and consider the following diagram
$$
\xymatrix{
\S_q^{\simeq} \h_n \C \ar[r]^{d_q} \ar[d] & \S_{q-1}^{\simeq} \h_n\C \ar[d] \\
\big((\h_n \C)^{\Delta^{q-1}}\big)^{\simeq} \ar[r]_{d_q}  & \big( (\h_n \C)^{\Delta^{q-2}}\big)^{\simeq}.
}
$$

For $0 \leq k \leq q$, let $T^q_k \subseteq \Ar[q]$ be the full subcategory which contains the subposet $\Ar[q-1] \subseteq \Ar[q]$ and the elements $\{(j \to q) | \ 0 \leq j \leq k \}$. Moreover, let $\mathscr{T}^q_k$ denote the full subcategory of $\Map(T^q_k,  \h_n \C)$ whose objects satisfy properties (i)--(ii) (Subsection \ref{K-theory_distinguished_squares}) restricted to $T^q_k$. In other words, this is the full subcategory (indeed an $\infty$-groupoid) which is spanned by the image of $\S_q^{\simeq} \h_n \C$ under the restriction map $\Map(\Ar[q], \h_n \C) \to \Map(T^q_k, \h_n \C)$. Then we may factorize the square above as the composition of the following squares
\begin{equation*} \label{diagram2}
\xymatrix{
\S_q^{\simeq} \h_n \C = \mathscr{T}^q_q \ar[r]^{\simeq} \ar[d] & \mathscr{T}^q_{q-1} \ar[r] \ar[d] & \cdots \ar[d] \ar[r] & \mathscr{T}^q_1 \ar[d] \ar[r] & \mathscr{T}^q_0 \ar[d]  \\ 
\big((\h_n \C)^{\Delta^{q-1}}\big)^{\simeq} \ar@{=}[r] & \big((\h_n \C)^{\Delta^{q-1}}\big)^{\simeq} \ar@{=}[r]  & \cdots \ar@{=}[r] & \cdots \ar@{=}[r] & \big((\h_n \C)^{\Delta^{q-1}}\big)^{\simeq}
}
\end{equation*}
followed by the pullback square
\begin{equation*} \label{diagram3}
\xymatrix{
\mathscr{T}^q_0 \ar[r] \ar[d] & \S_{q-1}^{\simeq} \h_n\C \ar[d] \\ 
\big((\h_n \C)^{\Delta^{q-1}}\big)^{\simeq} \ar[r]_{d_q} & \big( (\h_n \C)^{\Delta^{q-2}}\big)^{\simeq}.
}
\end{equation*}
We note that all the horizontal and vertical maps in these diagrams are given by the respective restriction functors. We claim that each map $\mathscr{T}^q_k \to \mathscr{T}^q_{k-1}$ is $n$-connected, for any $1 \leq k < q$, from which the required result follows. To see this, we consider the pullback of $\infty$-groupoids
\begin{equation} \label{diagram4}
\xymatrix{
\mathscr{T}^q_k \ar[d] \ar[r] & \big(\h_n(\C)^{\square, \can}\big)^{\simeq} \ar[d] \\
\mathscr{T}^q_{k-1} \ar[r] & \big(\h_n(\C)^{\ulcorner}\big)^{\simeq}
}
\end{equation} 
where the bottom map is given by the restriction to the subposet of $\mathscr{T}^q_{k-1}$
$$
\xymatrix{
(k-1 \to q-1) \ar[d] \ar[r] & (k-1 \to q) \\
(k \to q-1) 
}
$$ 
and $\big(\h_n(\C)^{\square, \can}\big)^{\simeq} \subset \big(\h_n(\C)^{\square}\big)^{\simeq}$ is the full subgroupoid that is spanned by the weak pushouts of order $n-1$. The right vertical map in \eqref{diagram4} is given by restriction along the upper corner inclusion in $\square = \Delta^1 \times \Delta^1$. The fiber of this map at $F \in \big(\h_n(\C)^{\ulcorner}\big)^{\simeq}$ is exactly the $\infty$-groupoid of weak colimits of $F$ 
of order $n-1$. Since $\h_n \C$ admits weak pushouts of order $n-1 > 0$, it follows that the fibers of the right vertical map in \eqref{diagram4} are 
$(n-1)$-connected. This means that the vertical maps in \eqref{diagram4} are $n$-connected and the required result follows.     
\end{proof}

\begin{example}
Theorem \ref{comparison_K-theory_ncat} for $n=1$ shows that the map $\rho_1 \colon K(\C) \to K(\h_1 \C, \can)$ is $1$-connected. In particular, this recovers the well-known fact that $K_0(\C)$ can be obtained from $\h_1 \C$ equipped with the canonical structure of distinguished squares.
\end{example} 

\begin{remark} \label{best_possible}
The connectivity estimate in Theorem \ref{comparison_K-theory_ncat} is best possible in general. Indeed, for $n=1$ and $\mathscr{E}$ an exact category, the comparison map $\rho_1$ for the $\infty$-category associated to the Waldhausen category of bounded chain complexes in $\mathscr{E}$ factors through the comparison map to Neeman's $K$-theory of the triangulated category $D^b(\mathscr{E})$, that is, we have maps
$$\rho_1 \colon K(\mathscr{E}) \xrightarrow{\beta \alpha} K({}^d D^b(\mathscr{E})) \to K(D^b(\mathscr{E}), \can)$$
where the last map is induced by the forgetful map of simplicial objects ${}^d \S_{\bullet} D^b(\mathscr{E}) \to \S_{\bullet} (D^b(\mathscr{E}), \can)$. 
(We refer to \cite{Neeman} for a nice overview of the $K$-theory of triangulated categories and details about the comparison maps $\alpha$ and  $\beta$.) The map induced on $K_1$ by $\beta \alpha$ is not injective in general by \cite[Section 11, Proposition 1]{Neeman}, see \cite[Sections 2 and 5]{Vaknin}.  
\end{remark}

Let us write $P_n X$ for the Postnikov $n$-section of a topological space $X$, which means that the canonical map $X \to P_n X$ is $(n+1)$-connected -- this agrees with the homotopy $n$-category of an $\infty$-groupoid. Then Theorem \ref{comparison_K-theory_ncat} implies that the functor $\C \mapsto P_{n-1}K(\C)$ descends to a functor defined for $(\h_n \C, \can)$. The following immediate corollary confirms a conjecture of Antieau in 
the case of connective $K$-theory \cite[Conjecture 8.35]{Antieau}.

\begin{corollary} \label{antieau_conj}
Let $\C$ and $\C'$ be stable $\infty$-categories such that there is an equivalence $(\h_n \C, \can) \simeq (\h_n \C', \can)$, as pointed $n$-categories with distinguished squares. Then there is a weak equivalence $$P_{n-1} K(\C) \simeq P_{n-1} K(\C').$$ 
\end{corollary} 

\begin{remark}
As explained in Remark \ref{canonical-structure-2}, the canonical structure on $(\h_n \C, \can)$ is preserved under equivalences for $n >1$, that is, an equivalence $\h_n \C \simeq \h_n \C'$ for $n>1$ 
is automatically an equivalence of pointed $n$-categories with distinguished squares. For $n=1$, note that an equivalence $\h_1 \C \stackrel{\Delta}{\simeq} \h_1 \C'$, as triangulated 
categories, clearly also respects the canonical structures of distinguished squares.
\end{remark}


\end{document}